\documentclass[letterpaper]{amsart} 
 
\usepackage{amssymb,latexsym} 
\usepackage{psfrag}
\usepackage{epsfig}
 
\numberwithin{equation}{section}

\newtheorem{theorem}{Theorem}[section] 
\newtheorem{proposition}[theorem]{Proposition}

\newtheorem{lemma}[theorem]{Lemma} 
 
\theoremstyle{definition}

\newtheorem{definition}[theorem]{Definition}

\begin{document}

\title{Cluster Algebras and Semipositive Symmetrizable Matrices}

\author{Ahmet I. Seven}

\address{Middle East Technical University, 06531, Ankara, Turkey}
\email{aseven@metu.edu.tr}

\thanks{The author's research was supported in part
by Turkish Scientific Research Council (TUBITAK)}

\date{November 18, 2009}


\begin{abstract}
There is a particular analogy between combinatorial aspects of cluster algebras and Kac-Moody algebras: roughly speaking, cluster algebras are associated with skew-symmetrizable matrices while Kac-Moody algebras correspond to (symmetrizable) generalized Cartan matrices. Both classes of algebras and the associated matrices have the same classification of finite type objects by the well-known Cartan-Killing types. In this paper, we study an extension of this correspondence to the affine type. In particular, we establish the cluster algebras which are determined by the generalized Cartan matrices of affine type.  
\end{abstract}


\maketitle

\section{Introduction}

\label{sec:introduction}

Cluster algebras are a class of commutative rings introduced by Fomin and Zelevinsky. It is well-known that these algebras are closely related with different areas of mathematics. A particular analogy exists between combinatorial aspects of cluster algebras and Kac-Moody algebras: roughly speaking, cluster algebras are associated with skew-symmetrizable matrices while Kac-Moody algebras correspond to (symmetrizable) generalized Cartan matrices. Both classes of algebras and the associated matrices have the same classification of finite type objects by the well-known Cartan-Killing types. In this paper, we study an extension of this correspondence between the two classes of matrices to the affine type. In particular, we establish the cluster algebras which are determined by the generalized Cartan matrices of affine type.  
 


To state our results, we need some terminology. In this paper, we deal with the combinatorial aspects of the theory of cluster algebras, so we will not need their definition nor their algebraic properties. The main combinatorial objects of our study will be skew-symmetrizable matrices and the corresponding directed graphs. Let us recall that an integer matrix $B$ is skew-symmetrizable if $DB$ is skew-symmetric for some diagonal matrix $D$ with positive diagonal entries. 
Recall also from \cite{CAII} that, for any matrix index $k$, the mutation of a skew-symmetrizable matrix $B$ in direction $k$ is another skew-symmetrizable matrix $\mu_k(B)=B'$ whose entries are given as follows: $B'_{i,j}=-B_{i,j}$ if $i=k$ or $j=k$; otherwise $B'_{i,j}=B_{i,j}+sgn(B_{i,k})[B_{i,k}B_{k,j}]_+$ (where we use the notation $[x]_+=max\{x,0\}$ and $sgn(x)=x/|x|$ with $sgn(0)=0$). Mutation is an involutive operation, so repeated mutations in all directions give rise to the \emph{mutation-equivalence} relation on skew-symmetrizable matrices. 
For each mutation (equivalence) class of skew-symmetrizable matrices, there is an associated cluster algebra \cite{CAII}. In this paper, we will establish the mutation-classes which are naturally determined by the generalized Cartan matrices of affine type. For this purpose, we use the following combinatorial construction from \cite{CAII}: for a skew-symmetrizable $n\times n$ matrix $B$, its \emph{diagram} is defined to be the directed graph $\Gamma (B)$ whose vertices are the indices $1,2,...,n$ such that there is a directed edge from $i$ to $j$ if and only if $B_{ij} > 0$, and this edge is assigned the weight $|B_{ij}B_{ji}|\,$. The diagram $\Gamma(B)$ does not determine $B$ as there could be several different skew-symmetrizable matrices whose diagrams are equal. In any case, we use the general term "diagram" to mean the diagram of a skew-symmetrizable matrix. Then the mutation $\mu_k$ can be viewed as a transformation on diagrams (see Section~\ref{sec:def} for a description).

On the other hand, an integer matrix $A$ is called symmetrizable if $DA$ is symmetric for some diagonal matrix $D$ with positive diagonal entries; we say that $A$ is (semi)positive if $DA$ is positive (semi)definite. Recall from \cite[Section~1]{BGZ} that a symmetrizable matrix $A$ is called a \emph{quasi-Cartan matrix} if all of its diagonal entries are equal to $2$. If the off-diagonal entries of a quasi-Cartan matrix are non-positive, then it is a generalized Cartan matrix; these are the matrices that give rise to Kac-Moody algebras, see \cite{K}. 
Motivated by the fact that cluster algebras and Kac-Moody algebras share the same classification of finite type objects \cite{CAII}, a notion of a \emph{quasi-Cartan companion} was introduced in \cite{BGZ} to relate skew-symmetrizable and symmetrizable matrices: a quasi-Cartan companion of a skew-symmetrizable matrix $B$ is a quasi-Cartan matrix $A$ such that $|A_{i,j}|= |B_{i,j}|$ for all $i \ne j$. In a slightly more general sense, we say that $A$ is a quasi-Cartan companion of a diagram $\Gamma$ if it is a quasi-Cartan companion of a skew-symmetrizable matrix whose diagram is equal to $\Gamma$. More combinatorially, a quasi-Cartan companion of a diagram may be viewed as a sign ($+$ or $-$) assignment to its edges (see Section~\ref{sec:def} for details). Given these definitions, it is natural to ask for an extension of the mutation operation on skew-symmetrizable matrices to their quasi-Cartan companions. One natural choice is the following \cite[Proposition~3.2]{BGZ}: for a  skew-symmetrizable matrix $B$ and a quasi-Cartan companion $A$, the "mutation of $A$ at $k$" is the quasi-Cartan matrix $A'$ such that, for any $i,j \ne k$, its entries are defined as  $A'_{i,k}=sgn(B_{i,k})A_{i,k}$, $A'_{k,j}=-sgn(B_{k,j})A_{k,j}$, $A'_{i,j}=A_{i,j}-sgn(A_{i,k}A_{k,j})[B_{i,k}B_{k,j}]_+$. It should be noticed that this definition uses both $B$ and $A$, so it can not be applied to an arbitrary quasi-Cartan matrix. Also the outcome $A'$, which is a quasi-Cartan matrix, may not be a quasi-Cartan companion of $\mu_k(B)=B'$.

In this paper, to identify a class of quasi-Cartan companions whose mutations are also quasi-Cartan companions, we introduce a notion of \emph{admissibility}. More specifically, for a skew-symmetrizable matrix $B$, we call a quasi-Cartan companion $A$ admissible if it satisfies the following sign condition: for any cycle $Z$ in $\Gamma(B)$, the product $\prod _{\{i,j\}\in Z}(-A_{i,j})$ over all edges of $Z$ is negative if $Z$ is oriented and positive if $Z$ is non-oriented; here a cycle\footnote{the term "chordless cycle" is used in \cite{BGZ}.} is an induced (full) subgraph isomorphic to a cycle (see Section~\ref{sec:def} for a precise definition). The main examples of admissible companions are the generalized Cartan matrices: if $\Gamma(B)$ is acyclic, i.e. has no oriented cycles at all, then the quasi-Cartan companion $A$ with $A_{i,j}= -|B_{i,j}|$, for all $i\ne j$, is admissible. However, for an arbitrary skew-symmetrizable matrix $B$, an admissible quasi-Cartan companion may not exist\footnote{it exists, e.g., if all cycles in $\Gamma(B)$ are cyclically oriented \cite[Corollary~5.2]{BGZ}.}. Our first result is a uniqueness property of these companions: if an admissible quasi-Cartan companion exists, then it is unique up to simultaneous sign changes in rows and columns (Theorem~\ref{th:adm unique}). 

To state our other results, we need to recall some more properties of companions. First let us note that the
admissibility property of a quasi-Cartan companion may not be preserved under mutations. However, it is preserved for some interesting classes of skew-symmetrizable matrices. The most basic class of such matrices are those of finite type. Recall from \cite{CAII} that a skew-symmetrizable matrix $B$ (or its diagram) is said to be of finite type if, for any $B'$ which is mutation-equivalent to $B$, we have $\left| B'_{i,j}B'_{j,i}\right| \leq 3$. Classification of finite type skew-symmetrizable matrices is identical to the famous Cartan-Killing classification \cite{CAII}. It follows that finite type skew-symmetrizable matrices can be characterized in terms of their diagrams as follows: $B$ is of finite type if and only if its diagram $\Gamma(B)$ is mutation-equivalent to a Dynkin diagram (Figure~\ref{fig:dynkin-diagrams}). Another characterization, which makes the relation to Cartan-Killing more explicit, was obtained in \cite{BGZ} using quasi-Cartan companions; in our setup it reads as follows: a skew-symmetrizable matrix $B$ is of finite type if and only if it has an admissible quasi-Cartan companion which is positive \cite[Theorem 1.2]{BGZ}. In particular, for a finite type skew-symmetrizable matrix, mutation of an admissible quasi-Cartan companion is also admissible. 

Given the finite type case, it is natural to ask for the relation between \emph{semipositive} symmetrizable matrices and skew-symmetrizable ones. In particular, it is natural to ask for an explicit description of the mutation classes of extended Dynkin diagrams (Figure~\ref{fig:extended-dynkin-diagrams}), which correspond to generalized Cartan matrices of affine type. In this paper we answer these questions and some others. We first show that each diagram in the mutation class of an extended Dynkin diagram has an admissible quasi-Cartan companion which is semipositive of corank $1$. However, unlike the finite type case, there exist other diagrams which have such a quasi-Cartan companion without being mutation-equivalent to any extended Dynkin diagram. We determine all those diagrams in Figure~\ref{fig:critical}; they appear in eight series depending on several parameters. In particular, we obtain the following description of the mutation class of an extended Dynkin diagram: a diagram $\Gamma$ is mutation-equivalent to an extended Dynkin diagram if and only if it has a semipositive admissible quasi-Cartan companion of corank $1$ and it does not contain any diagram which belongs to Figure~\ref{fig:critical} (Theorem~\ref{th:mut class ext}). We prove the theorem by showing that these two properties, when together, are invariant under mutations. In particular, we show that the mutation class of a skew-symmetrizable matrix $B$ whose diagram $\Gamma(B)$ is mutation-equivalent to an extended Dynkin diagram uniquely determines a generalized Cartan matrix of affine type (see Theorem~\ref{th:ext skew}).

After showing the existence of a semipositive admissible quasi-Cartan companion of corank $1$ on \emph{all} diagrams in the mutation class of an arbitrary extended Dynkin diagram, we show that the converse holds for diagrams of skew-symmetric matrices (i.e quivers\footnote{replacing an edge from a vertex $i$ to $j$ with weight $B_{i,j}^2$ by $B_{i,j}$ many arrows, the diagram of a skew-symmetric matrix $B$ can be viewed as a quiver}). More explicitly, we show that $\mathcal{S}$ is the mutation class of an extended Dynkin diagram corresponding to a skew-symmetric matrix if and only if every diagram in $\mathcal{S}$ has an admissible quasi-Cartan companion which is semipositive of corank $1$ (Theorem~\ref{th:mut class ext char}). 
Also we conjecture that any diagram in the mutation class of an acyclic diagram has an admissible quasi-Cartan companion which is equivalent to a generalized Cartan matrix (see Definition~\ref{def:q-equivalent} for the equivalence of quasi-Cartan matrices). 




In an important special case we prove stronger statements. To be more specific, let us first note that a semipositive quasi-Cartan companion of corank $1$ has a non-zero radical vector $u$; we call $u$ sincere if all of its coordinates are nonzero. We characterize the diagrams which have such a quasi-Cartan companion with a sincere radical vector as the diagrams of minimal infinite type (see Definition~\ref{def:min-2-infinite}, Theorem~\ref{th:minimal}). In particular, we show that these diagrams are  mutation-equivalent to an extended Dynkin diagram (see the theorem for a precise formulation). Diagrams of minimal infinite type were computed explicitly in \cite{S2} and their relation to cluster categories was studied in \cite{BRS}.


Given a diagram, one basic question is whether its mutation class is finite. It follows from our results that any extended Dynkin diagram has a finite mutation class. We also prove the converse: any acyclic diagram, with at least three vertices, which has a finite mutation class is either a Dynkin diagram or an extended Dynkin diagram (Theorem~\ref{th:acyclic fmc}). Thus we obtain another characterization of Dynkin and extended Dynkin diagrams. For diagrams of skew-symmetric matrices (i.e. quivers), this statement was obtained in \cite{BR} using categorical methods. In this paper we use more combinatorial methods for more general diagrams. Also the diagrams that we give in Figure~\ref{fig:critical} have finite mutation classes \cite{BS}; furthermore the diagrams from there which correspond to skew-symmetric matrices can be constructed from triangulations of surfaces as described in \cite{FST}. Thus it is natural to ask if the other diagrams in Figure~\ref{fig:critical} can be related to the approach in \cite{FST}.

The paper is organized as follows. In Section~\ref{sec:def}, we give basic definitions and prove our result on the uniqueness of an admissible quasi-Cartan companion for a diagram. In Section~\ref{sec:main-th}, we state our main results on the mutation classes of extended Dynkin diagrams and the associated quasi-Cartan companions. In Section~\ref{sec:pre}, we establish basic properties of (semipositive) admissible quasi-Cartan companions. In Section~\ref{sec:proof}, we prove our main results.












\section{Basic Definitions}

\label{sec:def}

In this section, we recall some definitions and statements from \cite{BGZ,BRS,CAII}. Throughout the paper, a matrix always means a square integer matrix.  

\begin{definition} 
\label{def:skew-symmetrizable} 
Let $B=(B_{i,j})$ be a $n \times n$ matrix (whose entries are integers). 
The matrix $B$ is called 
skew-symmetrizable if there exists a diagonal matrix $D$ 
with positive diagonal entries such that $DB$ is skew-symmetric. 
\end{definition}

\noindent
Skew-symmetrizable matrices can be characterized as follows \cite[Lemma~7.4]{CAII}: $B$ is skew-symmetrizable if and only if $B$ is sign-skew-symmetric (i.e. for any $i,j$ either $B_{i,j}=B_{j,i}=0$ or $B_{i,j}B_{j,i}<0$) and 
for all $k \geq 3$ and all 
$i_1, \dots, i_k\,$, it satisfies
\begin{equation} 
\label{eq:cycle=cycle}
B_{i_1,i_2} B_{i_2,i_3} \cdots B_{i_k,i_1} = 
(-1)^k B_{i_2,i_1} B_{i_3,i_2} \cdots B_{i_1,i_k}\,. 
\end{equation}
This characterization can be used conveniently in relation with the following construction which represents skew-symmetrizable matrices using graphs \cite[Definition~7.3]{CAII}:

\begin{definition} 
\label{def:diagram-of-B} 
Let $n$ be a positive integer and let $I=\{1,2,...,n\}$. 
The \emph{diagram} of a skew-symmetrizable (integer) 
matrix~$B=(B_{i,j})_{i,j\in I}$ is the weighted directed graph $\Gamma (B)$ 
with the vertex set $I$ such that there is a directed edge from $i$ to $j$ 
if and only if $B_{i,j} > 0$, and this edge is assigned the weight 
$|B_{i,j}B_{j,i}|\,$. 
\end{definition} 
\noindent
The property \eqref{eq:cycle=cycle} puts a condition on weights of graphs which represent skew-symmetrizable matrices.
To be more specific, let $\Gamma$ be as in the definition: a cycle $C$ in $\Gamma$ is an induced (full) subgraph whose vertices can be labeled by $\{1,2,...,r\},r\geq 3$, such that there is an edge between $i$ and $j$ if and only if $|i-j|=1$ or $\{i,j\}=\{1,r\}$. If the weights of the edges in $C$ are $w_1,w_2,...,w_r$, then the product $w_1w_2...w_r$ is a perfect square (i.e. square of an integer) by \eqref{eq:cycle=cycle}. 
Thus we can naturally define a \emph{diagram} as follows:

\begin{definition}
\label{lem:def-diagram}
A diagram $\Gamma$ is a finite directed graph (with no loops or $2$-cycles) whose edges are weighted with positive integers such that the product of weights along any cycle is a perfect square. 
\end{definition}
\noindent
By some abuse of notation, we denote by the same symbol $\Gamma$ 
the underlying undirected graph of a diagram. We denote an edge between vertices $i$ and $j$ by $\{i,j\}$. If $i$ is a vertex adjacent to an edge $e$, we sometimes say that "$i$ is on $e$". If an edge $e=\{i,j\}$ has weight which is equal to $1$, then we do not specify it in the picture. If all edges have weight $1$, then we call $\Gamma$ 
\emph{simply-laced}. 
By a \emph{subdiagram} of $\Gamma$, we always mean a diagram $\Gamma'$ 
obtained from $\Gamma$ by taking an induced (full) directed subgraph on a subset of
vertices and keeping all its edge weights the same as in $\Gamma$ 
\cite[Definition~9.1]{CAII}. We call a vertex $v$ \emph{source} (\emph{sink}) if all adjacent edges are oriented away (towards) $v$.
A diagram is called \emph{acyclic} if it has no oriented cycles at all. It is well-known that an acyclic diagram has a source and a sink.

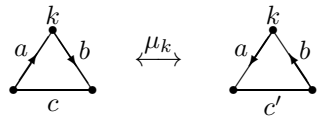
\begin{figure}[ht] 
\begin{center}
\setlength{\unitlength}{1.5pt} 
\begin{picture}(30,17)(-5,0) 
\put(0,0){\line(1,0){20}} 
\put(0,0){\line(2,3){10}} 
\put(0,0){\vector(2,3){6}} 
\put(10,15){\line(2,-3){10}} 
\put(10,15){\vector(2,-3){6}} 
\put(0,0){\circle*{2}} 
\put(20,0){\circle*{2}} 
\put(10,15){\circle*{2}} 
\put(2,10){\makebox(0,0){$a$}} 
\put(18,10){\makebox(0,0){$b$}} 
\put(10,-4){\makebox(0,0){$c$}} 
\put(10,19){\makebox(0,0){$k$}} 
\end{picture} 
$ 
\begin{array}{c} 
\stackrel{\textstyle\mu_k}{\longleftrightarrow} 
\\[.3in] 
\end{array} 
$ 
\setlength{\unitlength}{1.5pt} 
\begin{picture}(30,17)(-5,0) 
\put(0,0){\line(1,0){20}} 
\put(0,0){\line(2,3){10}} 
\put(10,15){\vector(-2,-3){6}} 
\put(10,15){\line(2,-3){10}} 
\put(20,0){\vector(-2,3){6}} 
\put(0,0){\circle*{2}} 
\put(20,0){\circle*{2}} 
\put(10,15){\circle*{2}} 
\put(2,10){\makebox(0,0){$a$}} 
\put(18,10){\makebox(0,0){$b$}} 
\put(10,-4){\makebox(0,0){$c'$}} 
\put(10,19){\makebox(0,0){$k$}} 
\end{picture} 
\end{center}
 
\vspace{-.2in} 
\caption{Diagram mutation} 
\label{fig:diagram-mutation-general} 
\end{figure}

For any vertex $k$ in a diagram $\Gamma$, there is the associated mutation $\mu_k$ which changes $\Gamma$ as follows:
\begin{itemize} 
\item The orientations of all edges incident to~$k$ are reversed, 
their weights intact. 
\item 
For any vertices $i$ and $j$ which are connected in 
$\Gamma$ via a two-edge oriented path going through~$k$ (see  
Figure~\ref{fig:diagram-mutation-general}), 
the direction of the edge $\{i,j\}$ in $\mu_k(\Gamma)$ and its weight $c'$ are uniquely determined by the rule 
\begin{equation} 
\label{eq:weight-relation-general} 
\pm\sqrt {c} \pm\sqrt {c'} = \sqrt {ab} \,, 
\end{equation} 
where the sign before $\sqrt {c}$ 
(resp., before $\sqrt {c'}$) 
is ``$+$'' if $i,j,k$ form an oriented cycle 
in~$\Gamma$ (resp., in~$\mu_k(\Gamma)$), and is ``$-$'' otherwise. 
Here either $c$ or $c'$ can be equal to~$0$, which means that the corresponding edge is absent. 
 
\item 
The rest of the edges and their weights in $\Gamma$ 
remain unchanged. 
\end{itemize} 

\noindent 
This operation is involutive, i.e. $\mu_k(\mu_k(\Gamma))=\Gamma$, so it defines an equivalence relation on the set of all diagrams. More precisely, two diagrams are called \emph{mutation-equivalent} if they can be obtained from each other by applying a sequence of mutations. The \emph{mutation class} of a diagram $\Gamma$ is the set of all diagrams which are mutation-equivalent to $\Gamma$. If $B$ is a skew-symmetrizable matrix, then $\Gamma(\mu_k(B))=\mu_k(\Gamma(B))$ (see 
Section~\ref{sec:introduction} for the definition of $\mu_k(B)$).

An important class of diagrams that behave very nicely under mutations are finite type diagrams:

\begin{definition} 
\label{def:2-finite} 
A diagram $\Gamma$ is said to be of \emph{finite type} if any diagram $\Gamma'$ 
which is mutation-equivalent to $\Gamma$ 
has all edge weights equal to $1,2$ or~$3$. A diagram is said to be of \emph{infinite type}
if it is not of finite type.
\end{definition}
\noindent
Let us note that a subdiagram of a finite type diagram is also of finite type. Every diagram which is mutation-equivalent to a diagram of finite type is of finite type itself. Also a diagram of finite type is of finite mutation type, i.e. its mutation class is finite. 

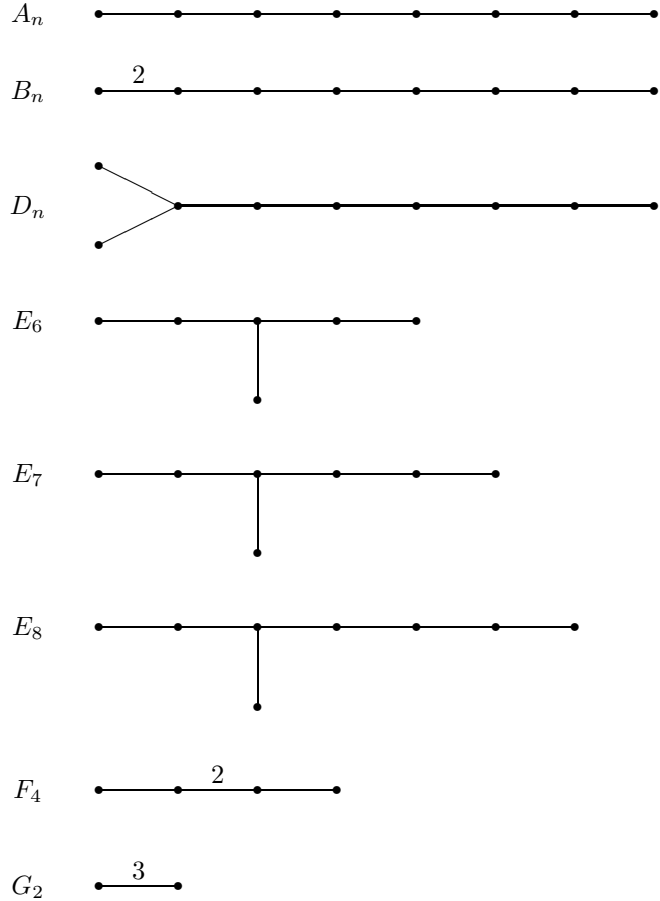
\begin{figure}[ht] 
\vspace{-.2in} 
\[ 
\begin{array}{ccl} 
A_n && 
\setlength{\unitlength}{1.5pt} 
\begin{picture}(140,17)(0,-2) 
\put(0,0){\line(1,0){140}} 
\multiput(0,0)(20,0){8}{\circle*{2}} 
\end{picture}\\
B_n
&& 
\setlength{\unitlength}{1.5pt} 
\begin{picture}(140,17)(0,-2) 
\put(0,0){\line(1,0){140}} 
\multiput(0,0)(20,0){8}{\circle*{2}} 
\put(10,4){\makebox(0,0){$2$}} 
\end{picture} 
\\[.2in] 
D_n 
&& 
\setlength{\unitlength}{1.5pt} 
\begin{picture}(140,17)(0,-2) 
\put(20,0){\line(1,0){120}} 
\put(0,10){\line(2,-1){20}} 
\put(0,-10){\line(2,1){20}} 
\multiput(20,0)(20,0){7}{\circle*{2}} 
\put(0,10){\circle*{2}} 
\put(0,-10){\circle*{2}} 
\end{picture} 
\\[.2in] 
E_6 
&& 
\setlength{\unitlength}{1.5pt} 
\begin{picture}(140,17)(0,-2) 
\put(0,0){\line(1,0){80}} 
\put(40,0){\line(0,-1){20}} 
\put(40,-20){\circle*{2}} 
\multiput(0,0)(20,0){5}{\circle*{2}} 
\end{picture} 
\\[.4in] 
E_7 
&& 
\setlength{\unitlength}{1.5pt} 
\begin{picture}(140,17)(0,-2) 
\put(0,0){\line(1,0){100}} 
\put(40,0){\line(0,-1){20}} 
\put(40,-20){\circle*{2}} 
\multiput(0,0)(20,0){6}{\circle*{2}} 
\end{picture} 
\\[.4in] 
E_8 
&& 
\setlength{\unitlength}{1.5pt} 
\begin{picture}(140,17)(0,-2) 
\put(0,0){\line(1,0){120}} 
\put(40,0){\line(0,-1){20}} 
\put(40,-20){\circle*{2}} 
\multiput(0,0)(20,0){7}{\circle*{2}} 
\end{picture} 
\\[.45in] 
F_4 
&& 
\setlength{\unitlength}{1.5pt} 
\begin{picture}(140,17)(0,-2) 
\put(0,0){\line(1,0){60}} 
\multiput(0,0)(20,0){4}{\circle*{2}} 
\put(30,4){\makebox(0,0){$2$}} 
\end{picture} 
\\[.1in] 
G_2 
&& 
\setlength{\unitlength}{1.5pt} 
\begin{picture}(140,17)(0,-2) 
\put(0,0){\line(1,0){20}} 
\multiput(0,0)(20,0){2}{\circle*{2}} 
\put(10,4){\makebox(0,0){$3$}} 
\end{picture} 
\end{array} 
\] 
\caption{Dynkin diagrams are arbitrary orientations of the Dynkin graphs given above; all orientations of the same Dynkin graph are mutation-equivalent to each other (this definition of a Dynkin diagram has been introduced in \cite{CAII}; note its difference from the definition in \cite{K}, where only the edges with multiple weights are oriented)} 
\label{fig:dynkin-diagrams} 
\end{figure}

Finite type diagrams were classified by Fomin and Zelevinsky in \cite{CAII}.
Their classification is identical to the Cartan-Killing classification. More precisely:

\begin{theorem}
\label{th:2-finite-class}
A connected diagram is of finite type if and only if it is mutation-equivalent to an arbitrarily oriented Dynkin diagram (Fig.~\ref{fig:dynkin-diagrams}). 
\end{theorem}

There is another description of finite type diagrams using the following notion:

\begin{definition} 
\label{def:min-2-infinite} 
A diagram $\Gamma$ is said to be of \emph{minimal infinite} type if it is of infinite type and any proper subdiagram of 
$\Gamma$ is of finite type.
\end{definition}

\noindent
A diagram is of finite type if and only if it does not contain any minimal infinite type diagram as a subdiagram. 
A complete list of minimal infinite type diagrams was obtained in \cite{S2}. We give a more algebraic characterization of these diagrams in Theorem~\ref{th:minimal}.

Another description of finite type diagrams was obtained in \cite{BGZ} using the following notion of "quasi-Cartan matrices", which we will use in this paper to describe the mutation classes of other types of diagrams:

\begin{definition} 
\label{def:symmetrizable} 
Let $A$ be a $n \times n$ matrix (whose entries are integers). 
The matrix $A$ is called symmetrizable if there exists a diagonal matrix $D$ with positive diagonal entries such that $DA$ is symmetric. We say that $A$ is a quasi-Cartan matrix if it is symmetrizable and all of its diagonal entries are equal to $2$.
\end{definition}
\noindent
The symmetrizable matrix $A$ is sign-symmetric, i.e. $sgn(A_{i,j})=sgn(A_{j,i})$. We say that $A$ is (semi)positive if $DA$ is positive (semi)definite, i.e. (resp. $x^TDAx\geq 0$) $x^TDAx>0$ for all $x\ne 0$ (here $x^T$ is the transpose of $x$ which is a vector viewed as a column matrix). We say that $u$ is a \emph{radical} vector of $A$ if $Au=0$; we call $u$ \emph{sincere} if all of its coordinates are non-zero. We call $A$ \emph{indefinite} if it is not semipositive. A quasi-Cartan matrix is a \emph{generalized Cartan matrix} if all of its non-zero entries which are not on the diagonal are negative.  

We use the following equivalence relation on quasi-Cartan matrices (recall that we work with matrices over integers):
\begin{definition} 
\label{def:q-equivalent} 
Quasi-Cartan matrices $A$ and $A'$ are called \emph{equivalent} if they have the same symmetrizer $D$, i.e $D$ is a diagonal matrix with positive diagonal entries such that both $C=DA$ and $C'=DA'$ are symmetric, and the symmetrized matrices satisfy $C'=E^TCE$ for some integer matrix $E$ with determinant $\mp 1$.
\end{definition}
\noindent
An important example of the equivalence for quasi-Cartan matrices is provided by the sign change operation: more specifically, the "\emph{sign change} at (vertex) $k$" replaces $A$ by $A'$ obtained by multiplying the $k$-th row and column of $A$ by $-1$. 

Quasi-Cartan matrices are related to skew-symmetrizable matrices via the following notion:
\begin{definition} 
\label{def:companion} 
Let $B$ be a skew-symmetrizable matrix. A \emph{quasi-Cartan companion} (or "companion" for short) of $B$ is a quasi-Cartan matrix $A$ with $|A_{i,j}|= |B_{i,j}|$ for all $i \ne j$. More generally, we say that $A$ is a quasi-Cartan companion of a diagram $\Gamma$ if it is a companion for a skew-symmetrizable matrix $B$ whose diagram is equal to $\Gamma$.

\end{definition}
\noindent
We define \emph{the restriction of the companion $A$ to a subdiagram} $\Gamma'$ as the quasi Cartan matrix obtained from $A$ by removing the rows and columns corresponding to the vertices which are not in $\Gamma'$. 
If $B$ is skew-symmetric, then any quasi-Cartan companion of it is symmetric; in this case we sometimes call $A_{i,j}$ the restriction of $A$ to the edge $\{i,j\}$. 

Let us note that for a diagram $\Gamma$, we may view a quasi-Cartan companion $A$ as a sign assignment to the edges (of the underlying undirected graph) of $\Gamma$; more explicitly an edge $\{i,j\}$ is assigned the sign of the entry $A_{i,j}$ (which is the same as the sign of $A_{j,i}$ because $A$ is sign-symmetric).


Motivated by the works in \cite{BGZ,BRS}, we introduce the following notion: 

\begin{definition} 
\label{def:admissible} 
Suppose that $B$ is a skew-symmetrizable matrix and let $A$ be a quasi-Cartan companion of $B$. We say that $A$ is \emph{admissible} if it satisfies the following sign condition: for any cycle $Z$ in $\Gamma$, the product $\prod _{\{i,j\}\in Z}(-A_{i,j})$ over all edges of $Z$ is negative if $Z$ is oriented and positive if $Z$ is non-oriented.   

\end{definition}
\noindent
The sign condition in the definition can also be described as follows: if $Z$ is (non)oriented, then there is exactly an (resp. even) odd number of edges $\{i,j\}$ such that $(A_{i,j})>0$. (recall that, since $A$ is symmetrizable, we have $sgn(A_{i,j})=sgn(A_{j,i})$). Thus an admissible quasi-Cartan companion distinguishes between the oriented and non-oriented cycles in a diagram. Note also that $A$ is admissible if and only if its restriction to any cycle is admissible. Thus the restriction of an admissible companion to a subdiagram is also admissible. Let us also note that sign change at  a vertex preserves admissibility.

In general, for a diagram $\Gamma$, an admissible quasi-Cartan companion may not exist. It is guaranteed to exist, e.g., if $\Gamma$ does not have any non-oriented cycles \cite[Corollary~5.2]{BGZ}. Our first result is that if an admissible companion exists, then it is  unique up to sign changes:


\begin{theorem}
\label{th:adm unique}
Suppose that $B$ is a skew-symmetrizable matrix. Let $A$ and $A'$ be any two admissible quasi-Cartan companions of $B$. Then $A$ and $A'$ can be obtained from each other by a sequence of simultaneous sign changes in rows and columns. In particular, $A$ and $A'$ are equivalent.
\end{theorem}
\noindent
This theorem generalizes \cite[Lemma~6.2]{BRS}. We will prove the theorem at the end of this section for convenience.

To proceed, let us first recall a characterization of finite type diagrams using quasi-Cartan companions, which reads in our setup as follows: 
\begin{theorem}
\label{th:2-finite-class positive}
\cite[Theorem~1.2]{BGZ} A diagram is of finite type if and only if it has an admissible quasi-Cartan companion which is positive.
\end{theorem}

\begin{figure}
\[ 
\begin{array}{ccl} 
\check{B_n}^{(4)} 
&& 
\setlength{\unitlength}{1.5pt} 
\begin{picture}(140,15)(-20,-5) 
\put(20,0){\line(1,0){120}} 
\put(0,10){\vector(2,-1){20}} 
\put(20,0){\vector(-2,-1){20}} 
\put(0,-10){\vector(0,1){20}} 
\put(-5,0){\makebox(0,0){$4$}}

\multiput(20,0)(20,0){7}{\circle*{2}} 
\put(0,10){\circle*{2}} 
\put(0,-10){\circle*{2}} 
\put(130,4){\makebox(0,0){$2$}}

\put(118,-5){\makebox(0,0){$a_m$}}
\put(24,-5){\makebox(0,0){$a_1$}}

\put(6,13){\makebox(0,0){$b_2$}}
\put(6,-13){\makebox(0,0){$b_1$}}

\put(145,-5){\makebox(0,0){$c_1$}}

\end{picture} 
\\[.1in] 
{ \check{B_n}^{(1)}(m,r)}
&& 
\setlength{\unitlength}{1.0pt} 
\begin{picture}(205,60)(-20,-5) 
\put(100,0){\circle*{2.0}}

\put(80,20){\circle*{2.0}}
\put(80,0){\circle*{2.0}} 
\put(120,0){\circle*{2.0}}
\put(140,0){\circle*{2.0}}
\put(160,0){\circle*{2.0}}
\put(180,0){\circle*{2.0}}
\put(200,0){\circle*{2.0}}
\put(80,20){\circle*{2.0}}
\put(60,40){\circle*{2.0}}
\put(40,40){\circle*{2.0}}
\put(20,20){\circle*{2.0}}
\put(40,-20){\circle*{2.0}}
\put(20,0){\circle*{2.0}}
\put(60,-20){\circle*{2.0}}

\put(80,0){\vector(1,0){20}}
\put(100,0){\line(1,0){100}}
\put(80,0){\line(-1,-1){20}}
\put(80,20){\vector(0,-1){20}}
\put(100,0){\vector(-1,1){20}}

\put(80,20){\line(-1,1){20}}
\put(60,40){\line(-1,0){20}}
\put(40,40){\line(-1,-1){20}}
\put(20,20){\line(0,-1){20}}
\put(20,0){\line(1,-1){20}}
\put(40,-20){\line(1,0){20}}

\put(190,-4){\makebox(0,0){2}}
\put(84,-6){\makebox(0,0){$b_1$}}
\put(66,43){\makebox(0,0){$b_3$}}
\put(16,-2){\makebox(0,0){$b_i$}}
\put(64,-25){\makebox(0,0){$b_r$}}
\put(85,24){\makebox(0,0){$b_2$}}
\put(-50,-12){\makebox(0,0){{ $m\geq 1,~r\geq 3$}}}

\put(-50,-27){\makebox(0,0){{$(m+r=n)$}}}

\put(103,-5){\makebox(0,0){$a_1$}}
\put(123,-5){\makebox(0,0){$a_2$}}
\put(177,-5){\makebox(0,0){$a_m$}}

\put(202,-7){\makebox(0,0){$c_1$}}

\end{picture} 
\\[.5in] 
{ \check{B}_n^{(1)}(r)}
&& 
\setlength{\unitlength}{1.0pt} 
\begin{picture}(200,50)(-20,-7) 
\put(80,0){\circle*{2.0}}

\put(100,0){\circle*{2.0}}
\put(80,0){\circle*{2.0}} 
\put(80,20){\circle*{2.0}}
\put(60,40){\circle*{2.0}}
\put(40,40){\circle*{2.0}}
\put(20,20){\circle*{2.0}}
\put(40,-20){\circle*{2.0}}
\put(20,0){\circle*{2.0}}
\put(60,-20){\circle*{2.0}}

\put(80,0){\vector(1,0){20}}
\put(80,0){\line(-1,-1){20}}
\put(80,20){\vector(0,-1){20}}
\put(100,0){\vector(-1,1){20}}
\put(80,20){\line(-1,1){20}}
\put(60,40){\line(-1,0){20}}
\put(40,40){\line(-1,-1){20}}
\put(20,20){\line(0,-1){20}}
\put(20,0){\line(1,-1){20}}
\put(40,-20){\line(1,0){20}}

\put(180,20){\makebox(0,0){$b_i$ is a source or sink}}
\put(180,7){\makebox(0,0){for some $3\leq i\leq r$}}

\put(90,-4){\makebox(0,0){$2$}}
\put(93,13){\makebox(0,0){$2$}}
\put(81,-7){\makebox(0,0){$b_1$}}
\put(66,42){\makebox(0,0){$b_3$}}
\put(16,-2){\makebox(0,0){$b_i$}}
\put(65,-25){\makebox(0,0){$b_r$}}
\put(83,25){\makebox(0,0){$b_2$}}
\put(106,0){\makebox(0,0){$c_1$}}
\put(-50,-14){\makebox(0,0){{$ r\geq 3$}}}


\end{picture}
\\[.4in]
\check{D}_n^{(4)} 
&& 
\setlength{\unitlength}{1.2pt} 
\begin{picture}(145,17)(-35,-2) 
\put(20,0){\line(1,0){100}} 
\put(20,0){\vector(-2,1){20}} 
\put(0,-10){\vector(2,1){20}} 
\put(120,0){\line(2,-1){20}} 
\put(120,0){\line(2,1){20}} 
\multiput(20,0)(20,0){6}{\circle*{2}} 
\put(0,10){\circle*{2}} 
\put(0,-10){\circle*{2}} 
\put(140,10){\circle*{2}} 
\put(140,-10){\circle*{2}} 

\put(0,10){\vector(0,-1){20}} 

\put(-5,2){\makebox(0,0){{$4$}}}

\put(118,-5){\makebox(0,0){$a_m$}}
\put(24,-5){\makebox(0,0){$a_1$}}

\put(6,15){\makebox(0,0){$b_2$}}
\put(6,-15){\makebox(0,0){$b_1$}}

\put(145,12){\makebox(0,0){$c_2$}}
\put(145,-12){\makebox(0,0){$c_1$}}

\end{picture} 
\\[.3in] 
{\check{D_n}^{(4)}(m,r)}
&& 
\setlength{\unitlength}{1.0pt}
\begin{picture}(180,43)(-20,-2)
\put(40,0){\circle*{2.0}}
\put(20,10){\circle*{2.0}}  \put(26,15){\makebox(0,0){$b_2$}}
\put(20,-10){\circle*{2.0}} \put(26,-15){\makebox(0,0){$b_1$}}
\put(60,0){\circle*{2.0}}
\put(80,0){\circle*{2.0}}
\put(100,0){\circle*{2.0}}
\put(120,0){\circle*{2.0}}
\put(140,0){\circle*{2.0}}
\put(160,20){\circle*{2.0}}
\put(160,0){\circle*{2.0}}
\put(180,40){\circle*{2.0}}
\put(200,40){\circle*{2.0}}
\put(220,20){\circle*{2.0}}
\put(220,0){\circle*{2.0}}
\put(180,-20){\circle*{2.0}}
\put(200,-20){\circle*{2.0}}

\put(158,23){\makebox(-5,2){$c_2$}}
\put(158,-3){\makebox(0,-3){$c_1$}}
\put(178,43){\makebox(-5,2){$c_3$}}
\put(178,-23){\makebox(0,-1){$c_r$}}
\put(-50,-10){\makebox(0,0){{ $m\geq 1,~r\geq 3$}}}

\put(-50,-25){\makebox(0,0){{ $(m+r=n-1)$}}}

\put(40,0){\line(1,0){100}}
\put(20,-10){\vector(2,1){20}}
\put(40,0){\vector(-2,1){20}}
\put(160,20){\vector(1,1){20}}
\put(160,0){\vector(0,1){20}}
\put(140,0){\vector(1,0){20}}
\put(160,20){\vector(-1,-1){20}}
\put(180,-20){\vector(-1,1){20}}

\put(160,20){\line(1,1){20}}

\put(180,40){\vector(1,0){20}}
\put(200,40){\vector(1,-1){20}}
\put(220,20){\vector(0,-1){20}}
\put(220,0){\vector(-1,-1){20}}
\put(200,-20){\vector(-1,0){20}}

\put(20,10){\vector(0,-1){20}}
\put(15,-1){\makebox(0,0){$4$}}

\put(138,-5){\makebox(0,0){$a_m$}}
\put(44,-5){\makebox(0,0){$a_1$}}

\end{picture}
\\[.2in] 
{ \check{D_n}^{(1)}(m,r)}
&& 
\setlength{\unitlength}{1.0pt} 
\begin{picture}(140,60)(-20,-2) 
\put(100,0){\circle*{2.0}}

\put(80,20){\circle*{2.0}}
\put(80,0){\circle*{2.0}} 
\put(120,0){\circle*{2.0}}
\put(140,0){\circle*{2.0}}
\put(160,0){\circle*{2.0}}
\put(180,0){\circle*{2.0}}
\put(80,20){\circle*{2.0}}
\put(60,40){\circle*{2.0}}
\put(40,40){\circle*{2.0}}
\put(20,20){\circle*{2.0}}
\put(40,-20){\circle*{2.0}}
\put(20,0){\circle*{2.0}}
\put(60,-20){\circle*{2.0}}

\put(100,0){\vector(-1,0){20}}
\put(80,0){\line(-1,-1){20}}
\put(80,0){\vector(0,1){20}}
\put(80,20){\vector(1,-1){20}}
\put(100,0){\line(1,0){80}}

\put(80,20){\line(-1,1){20}}
\put(60,40){\line(-1,0){20}}
\put(40,40){\line(-1,-1){20}}
\put(20,20){\line(0,-1){20}}
\put(20,0){\line(1,-1){20}}
\put(40,-20){\line(1,0){20}}

\put(180,0){\line(2,1){20}}
\put(180,0){\line(2,-1){20}}
\put(200,10){\circle*{2.0}}
\put(200,-10){\circle*{2.0}}
\put(205,-15){\makebox(0,0){{ $c_1$}}}
\put(205,15){\makebox(0,0){{ $c_2$}}}

\put(83,-6){\makebox(0,0){$b_1$}}
\put(66,43){\makebox(0,0){$b_3$}}
\put(16,-2){\makebox(0,0){$b_i$}}
\put(64,-25){\makebox(0,0){$b_r$}}
\put(85,24){\makebox(0,0){$b_2$}}
\put(-50,-10){\makebox(0,0){{ $m\geq 1,~r\geq 3$}}}

\put(-50,-25){\makebox(0,0){{ $(m+r=n-1)$}}}

\put(103,-5){\makebox(0,0){$a_1$}}
\put(123,-5){\makebox(0,0){$a_2$}}
\put(177,-5){\makebox(0,0){$a_m$}}

\end{picture} 
\\[.3in] 
{ \check{D}_n^{(1)}(m,r,s):}
&& 
\setlength{\unitlength}{1.0pt} 
\begin{picture}(224,54)(-20,-2) 

\put(-50,-10){\makebox(0,0){{ $m\geq 1,~r,s\geq 3$}}}
\put(-50,-25){\makebox(0,0){{ $(m+r+s=n+1)$}}}

\put(85,-5){\makebox(0,0){$b_1$}} 
\put(85,23){\makebox(0,0){$b_2$}} 
\put(15,24){\makebox(0,0){$b_i$}} 
\put(65,-23){\makebox(0,0){$b_r$}}

\put(158,-5){\makebox(0,0){$c_1$}} 
\put(156,24){\makebox(0,0){$c_2$}} 
\put(226,20){\makebox(0,0){$c_i$}} 
\put(178,-23){\makebox(0,0){$c_s$}} 

\put(103,-5){\makebox(0,0){$a_1$}} 
\put(140,-5){\makebox(0,0){$a_m$}}

\put(40,-20){\line(1,0){20}} 
\put(40,40){\line(1,0){20}} 
\put(20,0){\line(0,1){20}} 
\put(80,0){\vector(0,1){20}} 
\put(20,0){\line(1,-1){20}} 
\put(60,40){\line(1,-1){20}} 
\put(20,20){\line(1,1){20}} 
\put(60,-20){\line(1,1){20}} 
\put(80,20){\vector(1,-1){20}}
\put(100,0){\vector(-1,0){20}} 

\put(160,20){\vector(-1,-1){20}}
\put(140,0){\vector(1,0){20}}
\put(160,0){\vector(0,1){20}} 
\put(160,20){\vector(1,1){20}} 
\put(180,40){\vector(1,0){20}}
\put(200,40){\vector(1,-1){20}}
\put(220,20){\vector(0,-1){20}}
\put(220,0){\vector(-1,-1){20}}
\put(200,-20){\vector(-1,0){20}}
\put(180,-20){\vector(-1,1){20}}


\multiput(40,-20)(20,0){2}{\circle*{2}} 
\multiput(40,40)(20,0){2}{\circle*{2}} 
\multiput(20,0)(60,0){2}{\circle*{2}} 
\multiput(20,20)(60,0){2}{\circle*{2}} 


\put(110,0){\circle*{2}} 
\put(130,0){\circle*{2}} 
\put(160,20){\circle*{2}}
\put(160,0){\circle*{2}}
\put(180,-20){\circle*{2}}
\put(180,40){\circle*{2}} 
\put(200,-20){\circle*{2}}
\put(200,40){\circle*{2}} 
\put(220,0){\circle*{2}}
\put(220,20){\circle*{2}} 

\multiput(100,0)(20,0){3}{\circle*{2}} 
\put(80,0){\line(1,0){60}} 
\end{picture} 
\\[.16in] 
{ \check{D}_n^{(1)}(r),} 
&& 
\setlength{\unitlength}{1.0pt} 
\begin{picture}(200,60)(-20,-2) 
\put(-50,-10){\makebox(0,0){{ $r\geq 3$}}}

\put(40,-30){\line(1,0){20}} 
\put(40,30){\line(1,0){20}} 
\put(20,-10){\line(0,1){20}} 
\put(80,-10){\vector(0,1){20}} 
\put(20,-10){\line(1,-1){20}} 
\put(60,30){\line(1,-1){20}} 
\put(20,10){\line(1,1){20}} 
\put(60,-30){\line(1,1){20}} 
\put(80,10){\vector(1,-1){20}}

\put(20,10){\circle*{2}} 
\put(20,-10){\circle*{2}} 
\put(40,-30){\circle*{2}}
\put(40,30){\circle*{2}} 
\put(60,-30){\circle*{2}} 
\put(60,-10){\circle*{2}} 
\put(60,30){\circle*{2}} 
\put(80,-10){\circle*{2}} 
\put(80,10){\circle*{2}} 
\put(100,-10){\circle*{2}} 

\put(58,-14){\makebox(0,0){$a_1$}} 
\put(105,-14){\makebox(0,0){$c_1$}} 
\put(83,-16){\makebox(0,0){$b_1$}} 
\put(84,15){\makebox(0,0){$b_2$}} 
\put(15,15){\makebox(0,0){$b_i$}} 
\put(65,-34){\makebox(0,0){$b_r$}}

\put(100,-10){\vector(-1,0){20}} 
\put(60,-10){\vector(1,0){20}}  
\put(80,10){\vector(-1,-1){20}} 

\put(180,15){\makebox(0,0){$b_i$ is a source or sink}}
\put(180,2){\makebox(0,0){for some $3\leq i\leq r$}}

\end{picture}
\\[.4in] 

\end{array} 
\] 
\caption{Diagrams that do not appear in the mutation classes of extended Dynkin diagrams: undirected edges are assumed to be arbitrarily oriented with the condition that any cycle with an unspecified orientation is not cyclically oriented. (Each graph has $n+1$ vertices.)} 
\label{fig:critical} 
\end{figure}
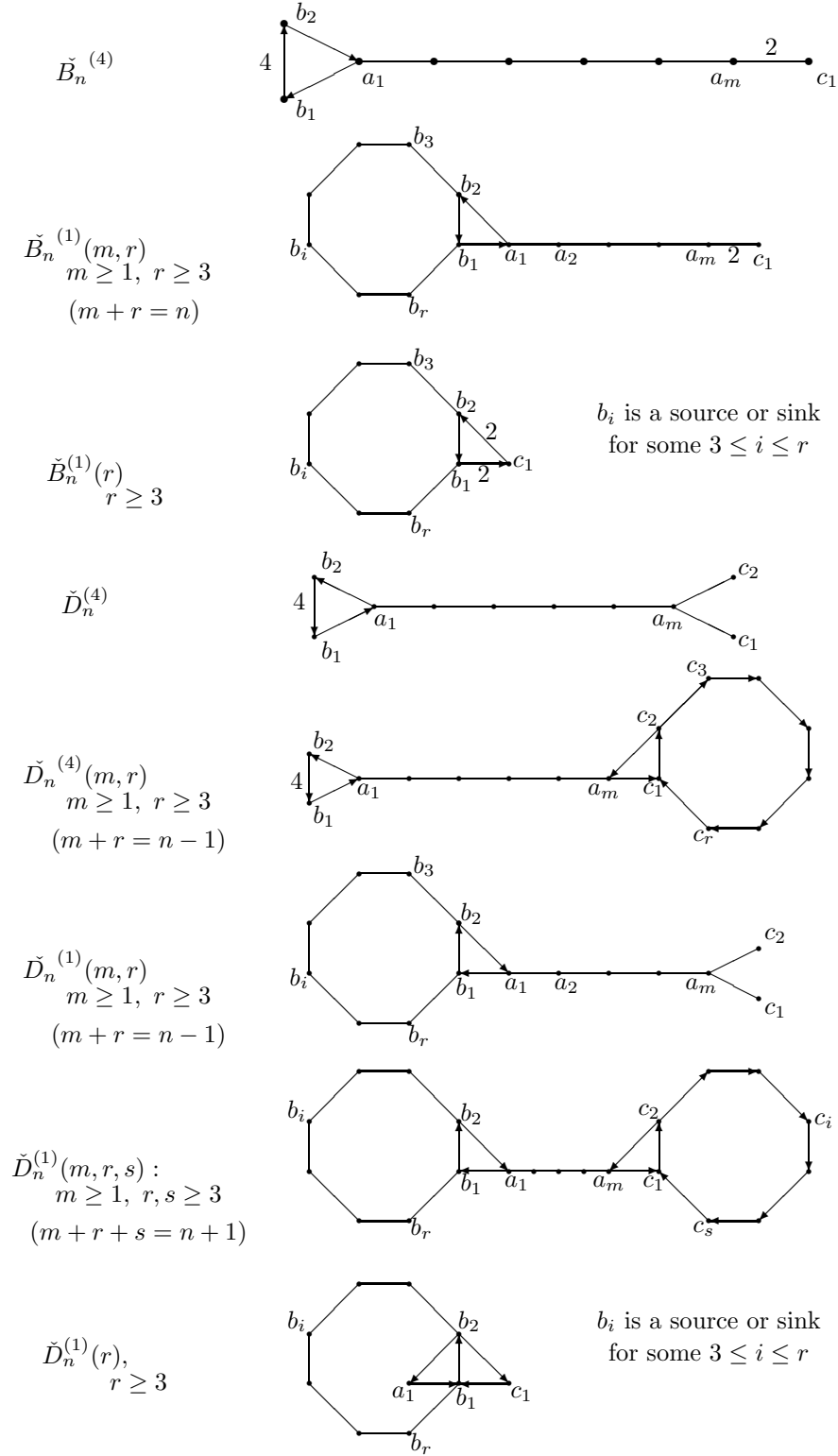

\noindent
The main tool in proving this theorem is the following operation on symmetrizable matrices analogous to the mutation operation on skew-symmetrizable matrices \cite[Proposition~3.2]{BGZ}:

\begin{definition} 
\label{def:comp-mut} 
Suppose that $\Gamma$ is a diagram and let $A$ be a quasi-Cartan companion of $\Gamma$. 
Let $k$ be a vertex in $\Gamma$. "The mutation of $A$ at $k$" is the quasi-Cartan matrix $A'$ such that for any $i,j \ne k$: $A'_{i,k}=sgn(B_{i,k})A_{i,k}$, $A'_{k,j}=-sgn(B_{k,j})A_{k,j}$, $A'_{i,j}=A_{i,j}-sgn(A_{i,k}A_{k,j})[B_{i,k}B_{k,j}]_+$. The quasi-Cartan matrix $A'$ is equivalent to $A$. It is a quasi-Cartan companion of $\mu_k(\Gamma)$ if $A$ is admissible \cite[Proposition~3.2]{BGZ}.
\end{definition}
\noindent
Note that $A'$ may not be admissible even if $A$ is admissible: e.g. if $A$ is an admissible quasi-Cartan companion of the diagram $\check{D}_5^{(4)}$ from Figure~\ref{fig:critical} and $k$ is the vertex $a_1$ there, then the corresponding $A'$ is not admissible. We conjecture that $A'$ is also admissible if $\Gamma$ is mutation-equivalent to an acyclic diagram (i.e. a diagram which has no oriented cycles at all). In this paper we prove this conjecture for the affine case, i.e. for diagrams which are mutation-equivalent to an extended Dynkin diagram. We will do more: we give an explicit description of their mutation classes and give some characterizing properties. 

\subsection{Proof of Theorem~\ref{th:adm unique}}
\label{sec:adm unique}

The theorem follows from the following lemma which is more general and stronger:

\begin{lemma}
\label{lem:adm unique}
Suppose that $B$ is a skew-symmetrizable matrix and let $\Gamma$ be the diagram of $B$. Let $A$ and $A'$ be any two (not necessarily admissible) quasi-Cartan companions of $B$. Suppose also that, for any cycle $C$ in $\Gamma$, the products  
$\prod _{\{i,j\}}(-A_{i,j})$ and $\prod _{\{i,j\}}(-A'_{i,j})$  over all edges of $C$ are equal. Then, viewing each $A$ and $A'$ as a sign assignment to the edges of $\Gamma$, we have the following: if $A$ and $A'$ are not equal, then  $A'$ can be obtained from $A$ by a sequence of sign changes at vertices such that a vertex is used at most once and not all vertices are used. 
\end{lemma}

We prove the lemma by induction on the number, say $n$, of vertices of $\Gamma$, which we can assume to be connected:  For $n=2$, the diagram $\Gamma$ has a single edge $e$; $A$ and $A'$ are not equal if they assign opposite signs to $e$, then sign change at any vertex transforms $A$ to $A'$.

Let us now assume that the lemma holds for diagrams with $n-1$ vertices or less. Let $\Delta$ be a connected subdiagram  obtained from $\Gamma$ by removing a vertex, say $n$ (the existence of such a vertex leaving a connected subdiagram is  easily seen). The vertices of $\Delta$ are $1,2,...,n-1$. 
Since $\Delta$ has less than $n$ vertices, by the induction argument we have the following: the restriction of $A$ to $\Delta$ can be transformed to the restriction of $A'$ using sign changes at vertices, say $1,...,r$, $r<n-1$ (i.e. as described in the lemma). Let $A''$ be the companion of $\Gamma$ obtained from $A$ by applying the same sign changes at $1,...,r$. Note that $A'_{i,j}=A''_{i,j}$ for all $i,j<n$ (i.e. $A'$ and $A''$ assign the same sign to any edge which is not adjacent to $n$). We claim that either $A'=A''$ or $A'$ can be obtained from $A''$ by a sign change at the vertex $n$. Note that, for any cycle $C$ in $\Gamma$, sign change at a vertex does not alter the product $\prod_{\{i,j\}}(-A_{i,j})$ over all edges of $C$, so $A'$ and $A''$ also satisfy the conditions of the lemma. If all edges ${\{i,n\}}$, $i<n$, are assigned the same sign by $A'$ and $A''$, then $A'=A''$ and we are done. If all edges ${\{i,n\}}$, $i<n$, are assigned opposite signs by $A'$ and $A''$, then $A'$ is obtained from $A''$ by the sign change at the vertex $n$, showing the lemma. The only remaining case then is the following: there are vertices $k$ and $m$ in $\Delta$, connected to the vertex $n$, such that the edge ${\{k,n\}}$ is assigned the same sign by both $A'$ and $A''$ but the edge ${\{m,n\}}$ is assigned opposite signs by them. Let us denote by $P$ a shortest path connecting $k$ and $m$ in $\Delta$. We can assume that $n$ is not connected to any vertex on $P$ other than $k$ and $m$ (otherwise we can find another pair of vertices like $k,m$ which are closer to each other). Then the subdiagram $\{P,n\}$ is a cycle such that exactly one of its edges is assigned opposite signs by $A'$ and $A''$; then, over the edges of this cycle, the products $\prod _{\{i,j\}}(-A'_{i,j})$ and $\prod _{\{i,j\}}(-A''_{i,j})$ are not equal, contradicting the assumption that $A'$ and $A''$ satisfy the condition of the lemma. This completes the proof.


\newpage

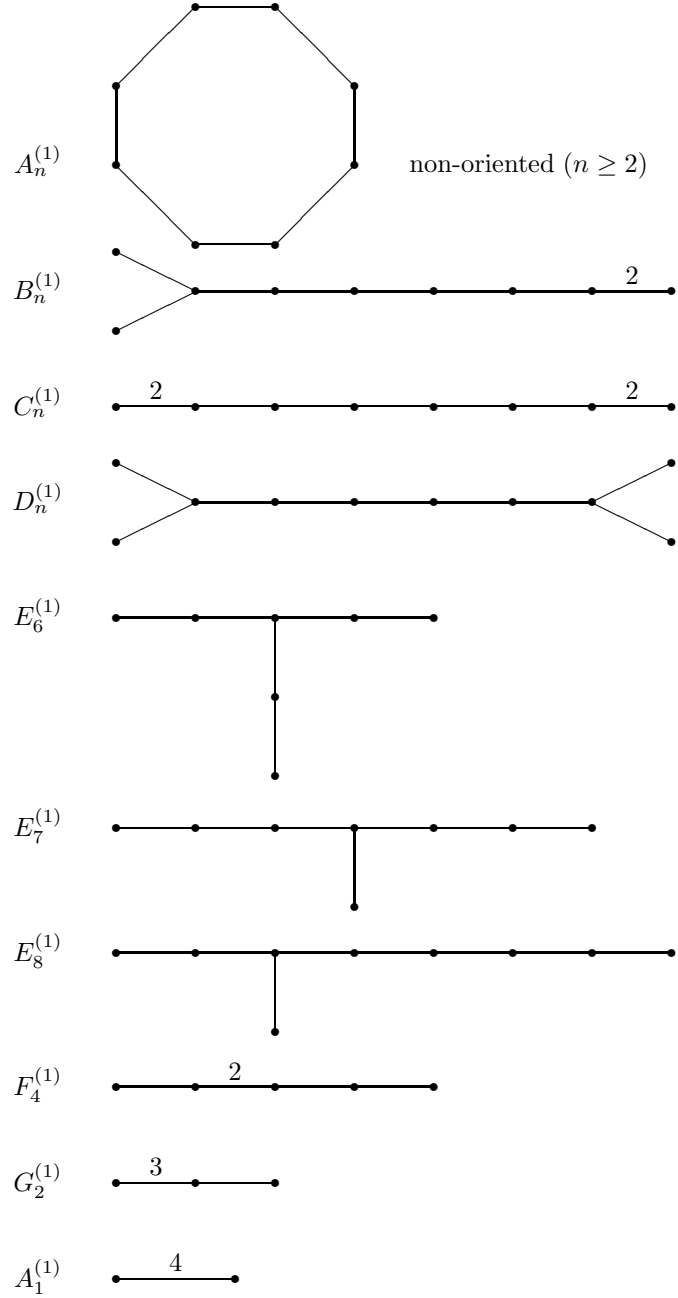
\begin{figure}[ht] 
\[ 
\begin{array}{ccl} 
A_n^{(1)}
&& 
\setlength{\unitlength}{1.5pt} 
\begin{picture}(140,60)(0,-2) 
\put(60,0){\circle*{2.0}} 
\put(60,20){\circle*{2.0}}
\put(40,40){\circle*{2.0}}
\put(20,40){\circle*{2.0}}
\put(0,20){\circle*{2.0}}
\put(20,-20){\circle*{2.0}}
\put(0,0){\circle*{2.0}}
\put(40,-20){\circle*{2.0}}

\put(60,0){\line(-1,-1){20}}
\put(60,0){\line(0,1){20}}
\put(60,20){\line(-1,1){20}}
\put(40,40){\line(-1,0){20}}
\put(20,40){\line(-1,-1){20}}
\put(0,20){\line(0,-1){20}}
\put(0,0){\line(1,-1){20}}
\put(20,-20){\line(1,0){20}}

\put(104,0){\makebox(0,0){non-oriented ($n\geq 2$)}}
\end{picture} 
\\[.3in] 
B_n^{(1)} 
&& 
\setlength{\unitlength}{1.5pt} 
\begin{picture}(140,15)(0,-2) 
\put(20,0){\line(1,0){120}} 
\put(0,10){\line(2,-1){20}} 
\put(0,-10){\line(2,1){20}} 
\multiput(20,0)(20,0){7}{\circle*{2}} 
\put(0,10){\circle*{2}} 
\put(0,-10){\circle*{2}} 
\put(130,4){\makebox(0,0){$2$}} 
\end{picture} 
\\[.2in] 
C_n^{(1)}
&& 
\setlength{\unitlength}{1.5pt} 
\begin{picture}(140,17)(0,-2) 
\put(0,0){\line(1,0){140}} 
\multiput(0,0)(20,0){8}{\circle*{2}} 
\put(10,4){\makebox(0,0){$2$}} 
\put(130,4){\makebox(0,0){$2$}} 
\end{picture} 
\\[.1in] 
D_n^{(1)} 
&& 
\setlength{\unitlength}{1.5pt} 
\begin{picture}(140,17)(0,-2) 
\put(20,0){\line(1,0){100}} 
\put(0,10){\line(2,-1){20}} 
\put(0,-10){\line(2,1){20}} 
\put(120,0){\line(2,-1){20}} 
\put(120,0){\line(2,1){20}} 
\multiput(20,0)(20,0){6}{\circle*{2}} 
\put(0,10){\circle*{2}} 
\put(0,-10){\circle*{2}} 
\put(140,10){\circle*{2}} 
\put(140,-10){\circle*{2}} 
\end{picture} 
\\[.2in] 
E_6^{(1)} 
&& 
\setlength{\unitlength}{1.5pt} 
\begin{picture}(140,17)(0,-2) 
\put(0,0){\line(1,0){80}} 
\put(40,0){\line(0,-1){40}} 
\put(40,-20){\circle*{2}} 
\put(40,-40){\circle*{2}} 
\multiput(0,0)(20,0){5}{\circle*{2}} 
\end{picture} 
\\[.7in] 
E_7^{(1)} 
&& 
\setlength{\unitlength}{1.5pt} 
\begin{picture}(140,17)(0,-2) 
\put(0,0){\line(1,0){120}} 
\put(60,0){\line(0,-1){20}} 
\put(60,-20){\circle*{2}} 
\multiput(0,0)(20,0){7}{\circle*{2}} 
\end{picture} 
\\[.25in] 
E_8^{(1)} 
&& 
\setlength{\unitlength}{1.5pt} 
\begin{picture}(140,17)(0,-2) 
\put(0,0){\line(1,0){140}} 
\put(40,0){\line(0,-1){20}} 
\put(40,-20){\circle*{2}} 
\multiput(0,0)(20,0){8}{\circle*{2}} 
\end{picture} 
\\[.3in] 
F_4^{(1)} 
&& 
\setlength{\unitlength}{1.5pt} 
\begin{picture}(140,17)(0,-2) 
\put(0,0){\line(1,0){80}} 
\multiput(0,0)(20,0){5}{\circle*{2}} 
\put(30,4){\makebox(0,0){$2$}} 
\end{picture} 
\\[.1in] 
G_2^{(1)} 
&& 
\setlength{\unitlength}{1.5pt} 
\begin{picture}(140,17)(0,-2) 
\put(0,0){\line(1,0){40}} 
\multiput(0,0)(20,0){3}{\circle*{2}} 
\put(10,4){\makebox(0,0){$3$}} 
\end{picture}\\[.1in] 
A_1^{(1)}
&& 
\setlength{\unitlength}{1.5pt} 
\begin{picture}(65,17)(20,-2) 
\put(20,0){\line(1,0){30}} 
\put(20,0){\circle*{2}} 
\put(50,0){\circle*{2}} 
\put(35,4){\makebox(0,0){$4$}} 
\end{picture} 
\end{array} 
\] 
\caption{Extended Dynkin diagrams are orientations of the extended Dynkin graphs given above; the first graph $A_n^{(1)}$ is assumed to be a  non-oriented cycle, the rest of the graphs are assumed to be arbitrarily oriented; each $X_n^{(1)}$ has $n+1$ vertices} 
\label{fig:extended-dynkin-diagrams} 
\end{figure} 

\clearpage

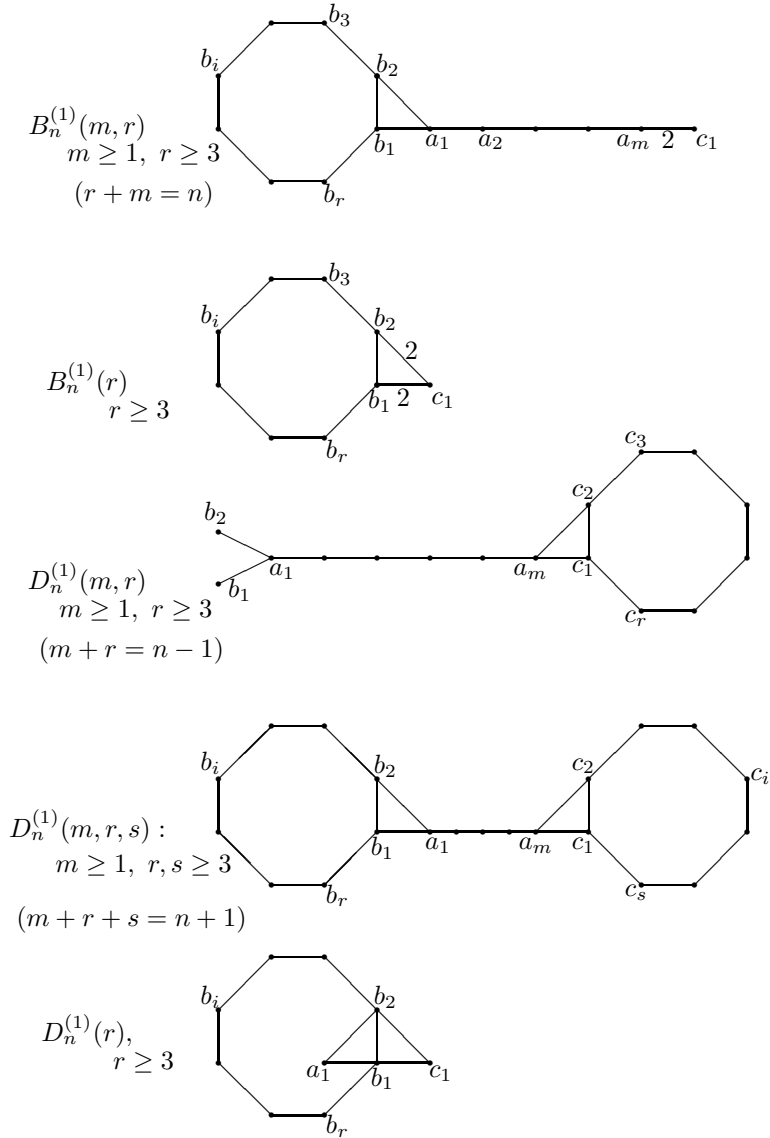
\begin{figure}
\[ 
\begin{array}{ccl} 
B_n^{(1)}(m,r)

&& 
\setlength{\unitlength}{1.0pt} 
\begin{picture}(205,60)(20,-2) 
\put(100,0){\circle*{2.0}}

\put(80,20){\circle*{2.0}}
\put(80,0){\circle*{2.0}} 
\put(120,0){\circle*{2.0}}
\put(140,0){\circle*{2.0}}
\put(160,0){\circle*{2.0}}
\put(180,0){\circle*{2.0}}
\put(200,0){\circle*{2.0}}
\put(80,20){\circle*{2.0}}
\put(60,40){\circle*{2.0}}
\put(40,40){\circle*{2.0}}
\put(20,20){\circle*{2.0}}
\put(40,-20){\circle*{2.0}}
\put(20,0){\circle*{2.0}}
\put(60,-20){\circle*{2.0}}

\put(80,0){\line(1,0){120}}
\put(80,0){\line(-1,-1){20}}
\put(80,0){\line(0,1){20}}
\put(100,0){\line(-1,1){20}}
\put(80,20){\line(-1,1){20}}
\put(60,40){\line(-1,0){20}}
\put(40,40){\line(-1,-1){20}}
\put(20,20){\line(0,-1){20}}
\put(20,0){\line(1,-1){20}}
\put(40,-20){\line(1,0){20}}

\put(190,-4){\makebox(0,0){2}}
\put(83,-6){\makebox(0,0){$b_1$}}
\put(66,43){\makebox(0,0){$b_3$}}
\put(17,26){\makebox(0,0){$b_i$}}
\put(64,-25){\makebox(0,0){$b_r$}}
\put(84,24){\makebox(0,0){$b_2$}}
\put(-10,-10){\makebox(0,0){{ $ m\geq 1,~r\geq 3$}}}

\put(-10,-25){\makebox(0,0){{ $(r+m=n)$}}}

\put(103,-5){\makebox(0,0){$a_1$}}
\put(123,-5){\makebox(0,0){$a_2$}}
\put(177,-5){\makebox(0,0){$a_m$}}
\put(205,-5){\makebox(0,0){$c_1$}}

\end{picture} 
\\[.6in] 
 B_n^{(1)}(r)
&& 
\setlength{\unitlength}{1.0pt} 
\begin{picture}(110,50)(20,-2) 
\put(80,0){\circle*{2.0}}

\put(100,0){\circle*{2.0}}
\put(80,0){\circle*{2.0}} 
\put(80,20){\circle*{2.0}}
\put(60,40){\circle*{2.0}}
\put(40,40){\circle*{2.0}}
\put(20,20){\circle*{2.0}}
\put(40,-20){\circle*{2.0}}
\put(20,0){\circle*{2.0}}
\put(60,-20){\circle*{2.0}}

\put(80,0){\line(1,0){20}}
\put(80,0){\line(-1,-1){20}}
\put(80,0){\line(0,1){20}}
\put(100,0){\line(-1,1){20}}
\put(80,20){\line(-1,1){20}}
\put(60,40){\line(-1,0){20}}
\put(40,40){\line(-1,-1){20}}
\put(20,20){\line(0,-1){20}}
\put(20,0){\line(1,-1){20}}
\put(40,-20){\line(1,0){20}}

\put(90,-5){\makebox(0,0){$2$}}
\put(93,13){\makebox(0,0){$2$}}
\put(81,-6){\makebox(0,0){$b_1$}}
\put(66,42){\makebox(0,0){$b_3$}}
\put(17,26){\makebox(0,0){$b_i$}}
\put(65,-25){\makebox(0,0){$b_r$}}
\put(83,25){\makebox(0,0){$b_2$}}
\put(-10,-10){\makebox(0,0){{$ r\geq 3$}}}
\put(105,-6){\makebox(0,0){$c_1$}}


\end{picture} 
\\[.4in] 
 D_n^{(1)}(m,r)
&& 
\setlength{\unitlength}{1.0pt}
\begin{picture}(180,43)(20,-12)
\put(40,0){\circle*{2.0}}
\put(20,10){\circle*{2.0}}  
\put(20,-10){\circle*{2.0}} 
\put(60,0){\circle*{2.0}}
\put(80,0){\circle*{2.0}}
\put(100,0){\circle*{2.0}}
\put(120,0){\circle*{2.0}}
\put(140,0){\circle*{2.0}}
\put(160,20){\circle*{2.0}}
\put(160,0){\circle*{2.0}}
\put(180,40){\circle*{2.0}}
\put(200,40){\circle*{2.0}}
\put(220,20){\circle*{2.0}}
\put(220,0){\circle*{2.0}}
\put(180,-20){\circle*{2.0}}
\put(200,-20){\circle*{2.0}}

\put(158,25){\makebox(0,0){$c_2$}}
\put(158,-4){\makebox(0,0){$c_1$}}
\put(178,45){\makebox(0,0){$c_3$}}
\put(178,-23){\makebox(0,0){$c_r$}}
\put(-13,-20){\makebox(0,0){{ $m\geq 1,~r\geq 3$}}}

\put(-13,-35){\makebox(0,0){{($m+r=n-1$)}}}

\put(40,0){\line(1,0){120}}
\put(40,0){\line(-2,-1){20}}
\put(40,0){\line(-2,1){20}}
\put(140,0){\line(1,1){20}}
\put(160,0){\line(0,1){20}}
\put(160,0){\line(1,-1){20}}
\put(160,20){\line(1,1){20}}
\put(180,40){\line(1,0){20}}
\put(200,40){\line(1,-1){20}}
\put(220,20){\line(0,-1){20}}
\put(220,0){\line(-1,-1){20}}
\put(200,-20){\line(-1,0){20}}

\put(138,-5){\makebox(0,0){$a_m$}}
\put(44,-5){\makebox(0,0){$a_1$}}

\put(28,-12){\makebox(0,0){$b_1$}}
\put(19,17){\makebox(0,0){$b_2$}}

\end{picture}
\\[.5in] 
{D_n^{(1)}(m,r,s):}
&& 
\setlength{\unitlength}{1.0pt} 
\begin{picture}(224,54)(20,-2) 

\put(-10,-13){\makebox(0,0){{ $m\geq 1,~r,s\geq 3$}}}

\put(-13,-33){\makebox(0,0){{($m+r+s=n+1$)}}}

\put(82,-6){\makebox(0,0){$b_1$}} 
\put(83,25){\makebox(0,0){$b_2$}} 
\put(17,26){\makebox(0,0){$b_i$}} 
\put(65,-24){\makebox(0,0){$b_r$}}

\put(158,-5){\makebox(0,0){$c_1$}} 
\put(158,25){\makebox(0,0){$c_2$}} 
\put(225,21){\makebox(0,0){$c_i$}} 
\put(178,-23){\makebox(0,0){$c_s$}} 

\put(103,-5){\makebox(0,0){$a_1$}} 
\put(140,-5){\makebox(0,0){$a_m$}}

\put(40,-20){\line(1,0){20}} 
\put(40,40){\line(1,0){20}} 
\put(20,0){\line(0,1){20}} 
\put(80,0){\line(0,1){20}} 
\put(20,0){\line(1,-1){20}} 
\put(60,40){\line(1,-1){20}} 
\put(20,20){\line(1,1){20}} 
\put(60,-20){\line(1,1){20}} 
\put(80,20){\line(1,-1){20}}
\put(140,0){\line(1,1){20}}
\put(140,0){\line(1,0){20}}
\put(160,0){\line(0,1){20}} 
\put(160,20){\line(1,1){20}} 
\put(180,40){\line(1,0){20}}
\put(200,40){\line(1,-1){20}}
\put(220,20){\line(0,-1){20}}
\put(220,0){\line(-1,-1){20}}
\put(200,-20){\line(-1,0){20}}
\put(180,-20){\line(-1,1){20}}

\put(40,-20){\line(1,0){12}} 
\put(60,40){\line(-1,0){12}} 
\put(20,20){\line(0,-1){12}} 
\put(80,0){\line(0,1){12}} 
\put(20,0){\line(1,-1){12}} 
\put(80,20){\line(-1,1){12}} 
\put(40,40){\line(-1,-1){12}} 
\put(60,-20){\line(1,1){12}}

\multiput(40,-20)(20,0){2}{\circle*{2}} 
\multiput(40,40)(20,0){2}{\circle*{2}} 
\multiput(20,0)(60,0){2}{\circle*{2}} 
\multiput(20,20)(60,0){2}{\circle*{2}} 


\put(110,0){\circle*{2}} 
\put(130,0){\circle*{2}} 
\put(160,20){\circle*{2}}
\put(160,0){\circle*{2}}
\put(180,-20){\circle*{2}}
\put(180,40){\circle*{2}} 
\put(200,-20){\circle*{2}}
\put(200,40){\circle*{2}} 
\put(220,0){\circle*{2}}
\put(220,20){\circle*{2}} 

\multiput(100,0)(20,0){3}{\circle*{2}} 
\put(80,0){\line(1,0){60}} 
\end{picture} 
\\[.19in] 
{ D_n^{(1)}(r),} 
&& 
\setlength{\unitlength}{1.0pt} 
\begin{picture}(140,60)(20,-2) 
\put(-10,-10){\makebox(0,0){{ $r\geq 3$}}}

\put(40,-30){\line(1,0){20}} 
\put(40,30){\line(1,0){20}} 
\put(20,-10){\line(0,1){20}} 
\put(80,-10){\line(0,1){20}} 
\put(20,-10){\line(1,-1){20}} 
\put(60,30){\line(1,-1){20}} 
\put(20,10){\line(1,1){20}} 
\put(60,-30){\line(1,1){20}} 
\put(80,10){\line(1,-1){20}}

\put(20,10){\circle*{2}} 
\put(20,-10){\circle*{2}} 
\put(40,-30){\circle*{2}}
\put(40,30){\circle*{2}} 
\put(60,-30){\circle*{2}} 
\put(60,-10){\circle*{2}} 
\put(60,30){\circle*{2}} 
\put(80,-10){\circle*{2}} 
\put(80,10){\circle*{2}} 
\put(100,-10){\circle*{2}} 

\put(58,-14){\makebox(0,0){$a_1$}} 
\put(104,-14){\makebox(0,0){$c_1$}} 
\put(82,-16){\makebox(0,0){$b_1$}} 
\put(83,15){\makebox(0,0){$b_2$}} 
\put(17,15){\makebox(0,0){$b_i$}} 
\put(65,-34){\makebox(0,0){$b_r$}}

\put(80,-10){\line(1,0){20}} 
\put(80,-10){\line(-1,0){20}}  
\put(60,-10){\line(1,1){20}} 
\end{picture}
\\[.6in]
\\[.1in] 

\end{array} 
\] 
\caption{Series of minimal infinite type diagrams which are not extended Dynkin: each graph above is assumed to have an arbitrary orientation such that all of its cycles are cyclically oriented (each $X_n^{(1)}$ has $n+1$ vertices)} 
\label{fig:minimal} 
\end{figure} 

\clearpage

\clearpage

\section{Main Results}
\label{sec:main-th}

We have already stated and proved one main result, Theorem~\ref{th:adm unique}, in Section~\ref{sec:def} for convenience.
In this section we state our remaining main results. We prove these results in Section~\ref{sec:proof} after some preparation in Section~\ref{sec:pre}.


Our first main result here is an explicit description of the mutation classes of extended Dynkin diagrams (Figure~\ref{fig:extended-dynkin-diagrams}):

\begin{theorem}
\label{th:mut class ext}
Let $B$ be a skew-symmetrizable matrix whose diagram $\Gamma(B)$ is connected. Then $\Gamma(B)$ is mutation-equivalent to an extended Dynkin diagram if and only if it does not contain any subdiagram that belongs to Figure~\ref{fig:critical} and $B$ has an admissible quasi-Cartan companion which is semipositive of corank $1$. 
\end{theorem}
\noindent
To use the theorem it is enough, by Theorem~\ref{th:adm unique}, to test just one admissible quasi-Cartan companion for semipositivity.



Our next result is the following classification statement as an analogue of \cite[Theorem 1.1]{BGZ}:

\begin{theorem}
\label{th:ext skew}
For a mutation class $\mathcal{S}$ of skew-symmetrizable matrices, the following are equivalent.
\begin{enumerate}
\item[(1)] There is a matrix in $\mathcal{S}$ whose diagram is mutation-equivalent to an extended Dynkin diagram.

\item[(2)] $\mathcal{S}$ contains a matrix $B$ with an admissible quasi-Cartan companion $A$ such that $A$ is a generalized Cartan matrix which is semipositive of corank $1$ (i.e. $A$ is of affine type \cite[Chapter 4]{K}). 
\end{enumerate}

Furthermore, the type of the generalized Cartan matrix in (2) is uniquely determined by $\mathcal{S}$.

Conversely, any generalized Cartan matrix of affine type except $A_n^{(1)}$, $n\geq 2$, uniquely determines a mutation class $\mathcal{S}$ of skew-symmetrizable matrices as in (1) (we refer to \cite[Chapter 4]{K} for a list of generalized Cartan matrices). 

\end{theorem}
\noindent
For skew-symmetric matrices, the second part this was obtained in \cite[Corollary~4]{CK} in a more general setup using cluster categories.

Our next result is the following characterization of extended Dynkin diagrams:

\begin{theorem}
\label{th:mut class ext char}
Let $\mathcal{S}$ be a mutation class of connected diagrams which correspond to skew-symmetric matrices. Then $\mathcal{S}$ is the mutation class of an extended Dynkin diagram if and only if every diagram in $\mathcal{S}$ has an admissible quasi-Cartan companion which is semipositive of corank $1$. 
\end{theorem}

\noindent
Let us note that this statement may be viewed as a converse of Theorem~\ref{th:mut class ext} for diagrams of skew-symmetric matrices (i.e quivers); it may not be true for diagrams of non-skew-symmetric matrices as can be checked on diagrams from Figure~\ref{fig:critical}. The crucial component in both theorems is the admissibility property, which is not preserved under mutation in general (Definition~\ref{def:comp-mut}) 
, however it is preserved in the situation of the theorems. More generally, we conjecture that admissibility property is preserved in the mutation class of any acyclic diagram.


Let us recall that a semipositive quasi-Cartan companion of corank $1$ has a non-zero radical vector $u$; we call $u$ \emph{sincere} if all of its coordinates are nonzero. We characterize all diagrams which have such a quasi-Cartan companion as follows:



\begin{theorem}
\label{th:minimal}

Let $\Gamma$ be a diagram with at least five vertices. Then $\Gamma$ is of minimal infinite type if and only if it has an admissible quasi-Cartan companion which is semipositive of corank $1$ with a sincere radical vector. Furthermore, if $\Gamma$ is of minimal infinite type, then it is mutation-equivalent to an extended Dynkin diagram. (If $\Gamma$ corresponds to a skew-symmetric matrix, then it is enough to have three vertices for the statements to be true).

\end{theorem}


Given a diagram, one basic question is whether its mutation class is finite. We determine all acyclic diagrams whose mutation classes are finite:   


\begin{theorem}
\label{th:acyclic fmc}
Let $\Gamma$ be an acyclic connected diagram with at least three vertices. Then the mutation class of $\Gamma$ is finite if and only if $\Gamma$ is either a Dynkin diagram or an extended Dynkin diagram. 
\end{theorem}
\noindent
For diagrams of skew-symmetric matrices (i.e. quivers), this statement was obtained in \cite{BR} using categorical methods. In this paper we use more combinatorial methods for more general diagrams. Let us also mention that there are algorithms to check whether a given skew-symmetric matrix is of finite mutation type: one of them is realized in B. Keller's computer program (which is available at www.math.jussieu.fr/\verb ~ keller/quivermutation); a polynomial-time algorithm is given in \cite{FSTu}.

\section{Preliminary results}
\label{sec:pre}

In this section we give some properties of semipositive quasi-Cartan companions. Their most basic property that we will use is the following:
\begin{proposition}\label{prop:semi-sub}
Suppose that $A$ is a semipositive quasi-Cartan companion of a diagram $\Gamma$.
Suppose also that $u$ is a radical vector for the restriction of $A$ to a subdiagram $\Sigma$, i.e.
$u$ is in the span of the standard basis vectors which correspond to the vertices in $\Sigma$ and 
$x^TAu=0$ for all $x$ in the same span. Then $u$ is a radical vector for $A$ as well (i.e. $x^TAu=0$ for \emph{all} $x$).   
\end{proposition}
\noindent
This statement is well-known. However, we could not find a suitable reference, therefore we give a proof here: let us assume that $u$ is not a radical vector for $A$. We can assume, without loss of generality, that $u\ne 0$. Then there is a vertex $k$ such that $e_k^TAu\ne 0$ (here $e_k$ is the $k$-th standard basis vector). Let $D=diag(d_1,...,d_n)$, $d_i>0$ for all $i$, be the symmetrizing matrix for $A$. We have $e_k^TDAu\ne 0$ as well (because $e_k^TD=d_ke_k^T$); assume without loss of generality that this number is negative (otherwise take $-e_k$ instead of $e_k$), note then that $e_k^TDAu\leq -1$ because we work over integers. Also note that, since $DA$ is symmetric, we have $e_k^TDAu=u^TDAe_k$. Let $a=e_k^TDAe_k$, which is positive (because it is equal to $d_kA_{k,k}=2d_k$). Then, e.g., for the vector $w=au+e_k$ we have $w^TDAw< 0$, contradicting that $A$ is semipositive. This completes the proof. 

Let us give some other properties of semipositive quasi-Cartan companions.
\begin{proposition}\label{prop:semidefinite comp}
Let $\Gamma$ be a diagram. 
Suppose that $A$ is a quasi-Cartan companion of $\Gamma$ which is semipositive. Then we have the following:

\begin{enumerate}
\item[(i)] The weight of any edge is at most $4$.

\item[(ii)] The restriction of $A$ to any edge of weight $4$ is not positive.

\item[(iii)] If $e$ is any edge whose weight is $4$, then any three-vertex diagram that contains $e$ is a triangle whose  edge weights are either $4,1,1$ or $4,4,4$ or $4,2,2$ or $4,3,3$.

\item[(iv)] If $C$ is a non-simply-laced cycle, then the product $\prod_{\{i,j\}\in C} (-A_{i,j})$ over all edges of $C$ is negative (so an odd number of edges of $C$ are assigned $(+)$ by $A$). 

\item[(v)] Suppose that $C$ is a simply-laced cycle such that for each edge the corresponding entry of $A$ is $-1$. Let $u$ be the vector whose coordinates are $1$ in the vertices of $C$ and $0$ in the remaining vertices. 
Then $u$ is a radical vector for $A$.

\item[(vi)] Suppose that $C$ is a simply-laced cycle such that the product $\prod_{\{i,j\}\in C} (-A_{i,j})$ over all edges of $C$ is positive. If a vertex $k$ is connected to $C$, then it is connected to at least two vertices in $C$.

\item[(vii)] Suppose that $\Gamma$ is {simply-laced} and let $C$ be a cycle in $\Gamma$ such that the product $\prod_{\{i,j\}\in C} (-A_{i,j})$ over all edges of $C$ is positive. If a vertex is connected to $C$, then 
it is connected to exactly an even number of vertices in $C$. 


\end{enumerate}

\end{proposition}

\noindent
Statements (i)-(v) easily follow from the definitions and known facts on generalized Cartan matrices \cite[Chapter~4]{K}. 
For (vi): applying sign changes if necessary (Theorem~\ref{th:adm unique}), we can assume that $C$ is as in part (v)  with the radical vector $u$. However, if $k$ connected to exactly one vertex in $C$ then $e_k^TAu\ne 0$, contradicting (v) (here $e_k$ is the $k$-th standard basis vector). Part (vii) is also proved similarly: assuming $C,u$ as in part (v), if $k$ is connected to exactly an odd number of vertices, then, for the edges connecting $k$ to $C$, the number of such edges assigned ($+$) is different from those assigned ($-$), implying that $e_k^TAu\ne 0$, which contradicts (v).

Let us now give some properties of admissible quasi-Cartan companions:

\begin{proposition}\label{prop:double in admissible}
Let $\Gamma$ be a diagram. Suppose that $A$ is an admissible quasi-Cartan companion which is semipositive. Then we have the following:
\begin{enumerate}
\item[(i)] If  $e$ is an edge whose weight is $4$, then any three-vertex subdiagram that contains $e$ is an \emph{oriented} triangle (see also part (iii) in the above proposition).

\item[(ii)] Any non-oriented cycle $C$ is simply-laced. Furthermore, the restriction of $A$ to $C$ is not positive.

\item[(iii)] Suppose that $A$ is of corank $1$ and let $i$ be a vertex which is on an edge whose weight is $4$ or on a non-oriented cycle. Then the subdiagram obtained by removing $i$ is of finite type.

\item[(iv)] Any diagram in Figures~\ref{fig:extended-dynkin-diagrams},~\ref{fig:minimal} 
has an admissible  quasi-Cartan companion of corank $1$ with a sincere radical vector.

\item[(v)] Suppose that $A$ is of corank $1$. Then $\Gamma$ contains at most one diagram from 
Figure~\ref{fig:extended-dynkin-diagrams} or Figure~\ref{fig:minimal} as a subdiagram. This is true, in particular, if $\Gamma$ contains an edge whose weight is $4$ or contains a non-oriented cycle.
\end{enumerate}
\end{proposition}

\noindent
These statements also follow easily from the definitions and known facts on generalized Cartan matrices \cite[Chapter~4]{K}. The  admissible quasi-Cartan companions of the diagrams in (iv) have also been studied in \cite{BRS}. Statement (v) follows from  Proposition~\ref{prop:semi-sub} and part (iv). 

Let us now look into the mutation operation given in Definition~\ref{def:comp-mut}. Recall that mutation of an admissible quasi-Cartan companion is also a quasi-Cartan companion, however it is not necessarily admissible. Our next statement gives one case when it is guaranteed to be admissible:


\begin{proposition}
\label{prop:mut adm}
Let $\Gamma$ be a diagram which does \emph{not} have any non-oriented cycles nor any edge whose weight is greater than or equal to $4$. Suppose that $A$ is an admissible quasi-Cartan companion of $\Gamma$ and let $A'$ be the quasi-Cartan  companion for $\mu_k(\Gamma)=\Gamma'$ obtained by mutating $A$ as in Definition~\ref{def:comp-mut}. Then $A'$ is also admissible.
\end{proposition}

\noindent
To prove this statement, we will need the following two lemmas which can be checked easily using the definitions:

\begin{lemma}\label{lem:C e mut}
Suppose that $\Gamma$ is a diagram which has at least three vertices and let $k$ be a vertex of $\Gamma$.
If $k$ is on a non-oriented cycle or on an edge whose weight is greater than or equal to $4$, then $\mu_k(\Gamma)$ contains an edge whose weight is at least $4$ or contains a non-oriented cycle.


\end{lemma}

\begin{lemma}\label{lem:C positive}
Let $C$ be a cycle (oriented or not). Let $Ck$ be 
a diagram obtained by connecting a new vertex $k$ to $C$ and let $A$ be a companion of $Ck$ such that 
the product $\prod_{\{i,j\}\in C} (-A_{i,j})$ is negative.
Suppose that $k$ is connected to an even number of vertices in $C$. Suppose also that
$k$ is connected to $C$ in such a way that it is connected to two vertices
which are not connected to each other in $C$ (this condition excludes only the case
when $k$ is connected to exactly two vertices in $C$ and those vertices are
connected to each other). Then $Ck$ necessarily has a 
cycle $C'$ which contains $k$ such that $\prod_{\{i,j\}\in C'} (-A_{i,j})$ is positive.
\end{lemma}

\noindent
{\bf Proof of Proposition~\ref{prop:mut adm}.} 
Let us denote by $A''$ the companion obtained by mutating $A'$ at $k$. Then $A''$ is a companion of $\mu_k(\Gamma')=\Gamma$ which is equal to $A$ up to a sign change at $k$. In particular, $A''$ is admissible. To prove the proposition, it is enough to show the following statement:

(***) if $A'$ is not admissible, i.e. there is a cycle $Z$ which does not satisfy the sign condition in Definition~\ref{def:admissible}, then $A''$ is not an admissible companion of $\mu_k(\Gamma')=\Gamma$ or $\Gamma$ contains a subdiagram which is a double edge or a non-oriented cycle .

To show (***), we first consider the case where $k$ is on $Z$. Note that if $k$ is a source or sink of $Z$, then $\mu_k(Z)$ is also a cycle on which $A''$ does not satisfy the same condition of admissibility. If $k$ is not a source or sink, then either $A''$ is not a companion of $\mu_k(\Gamma')$ (this happens when $Z$ is a triangle) or the diagram obtained from $\mu_k(Z)$ by removing $k$ is a cycle such that the restriction of $A''$ on it is not admissible, so $A''$ is not admissible.

We proceed by considering $k$ which is not on $Z$. Note that, by Lemma~\ref{lem:C e mut}, we can assume that any edge that is adjacent to $k$ has weight less than $4$ and any cycle $C'$ that contains $k$ is oriented and, by what we have have considered above, the restriction of $A'$ on it is admissible. For convenience, we will denote the subdiagram $\{Z,k\}$ by $Zk$.


{Case 1.} \emph{$Z$ is an oriented cycle.} 
If $k$ is connected to exactly one vertex in $Z$, then $\mu_k$ does not affect $Z$. Also if $k$ is connected to two vertices in $Z$ which are not connected to each other, then there is necessarily a non-oriented cycle that contains $k$ (because $Z$ is oriented), contradicting our assumption that any cycle that contains $k$ is oriented. Thus, for the rest of this case, we assume that $k$ is connected to exactly two vertices $z_1,z_2$ in $Z$ and $z_1,z_2$ are connected. By our assumption that any cycle that contains $k$ is oriented, the triangle $\{k,z_1,z_2\}$ is oriented. Let $w$ be the weight of the edge $\{z_1,z_2\}$ and let $p$ be the product of the weights of the edges $\{k,z_1\}$ and $\{k,z_2\}$. Then we have the following: if $p<w$, then in $\mu_k(\Gamma')=\Gamma$ the subdiagram $\{z_1,z_2,k\}$ is a non-oriented triangle; if $p=w$, then $\mu_k$ destroys the edge $\{z_1,z_2\}$, so $\mu_k(Zk)\subset \Gamma$ is an oriented cycle such that the restriction of $A''$ on it is not admissible; if $p>w$, then $\mu_k$ reverses the edge $\{z_1,z_2\}$, so in $\Gamma$ the subdiagram on $Z$ is a non-oriented cycle; in each case (***) holds.

{Case 2.} \emph{$Z$ is a non-oriented cycle.}  
If $k$ is connected to exactly one vertex in $Z$, then $\mu_k$ does not affect $Z$. Also if $k$ is connected to exactly an odd number $\geq 3$ vertices in $Z$, then there is necessarily a non-oriented cycle that contains $k$, contradicting our assumption that any cycle that contains $k$ is oriented.
Thus for the rest of this case we can assume that $k$ is connected to exactly an even number of vertices in $Z$. 
If $k$ is connected to two vertices in $Z$ which are not connected to each other, then by Lemma~\ref{lem:C positive} there is necessarily a cycle $C'$ that contains $k$ such that $\prod_{\{i,j\}\in C'} (-A'_{i,j})$ is positive, so $C'$ is non-oriented (because we assumed that the restriction of $A'$ to any cycle that contains $k$ is admissible), which contradicts our assumption that any cycle that contains $k$ is oriented. It remains to consider the subcase where $k$ is connected to exactly two vertices, say $z_1,z_2$, and $z_1,z_2$ are connected. By our assumption that any cycle that contains $k$ is oriented, the triangle $\{k,z_1,z_2\}$ is oriented. As in Case 1 above, let $w$ be the weight of the edge $\{z_1,z_2\}$ and let $p$ be the product of the weights of the edges $\{k,z_1\}$ and $\{k,z_2\}$. Then we have the following: if $p<w$, then in $\mu_k(\Gamma')=\Gamma$ the subdiagram $\{z_1,z_2,k\}$ is a non-oriented triangle; if $p=w$, then $\mu_k$ destroys the edge $\{z_1,z_2\}$, so $\mu_k(Zk)\subset \Gamma$ is a non-oriented cycle; if $p>w$, then $\mu_k$ reverses the edge $\{z_1,z_2\}$, so in $\Gamma$ either the subdiagram on $Z$ is a non-oriented cycle (this happens if there is a vertex $v\ne z_1,z_2$ such that $v$ is a source or sink in $Z$) or it is an oriented cycle 
such that the restriction of $A''$ on it is not admissible; in each case (***) holds. This completes the proof of Proposition~\ref{prop:mut adm}.


\section{Proofs of Main Results}
\label{sec:proof}

\subsection{Proof of Theorem~\ref{th:mut class ext}}
\label{sec:mut class ext}

For convenience we first prove the following statement:
\begin{proposition}
\label{prop:mut class ext}
Suppose that $\Gamma$ is a diagram which does not contain any subdiagram that belongs to Figure~\ref{fig:critical}. 
Let $A$ be an admissible quasi-Cartan companion of $\Gamma$ which is semipositive of corank $1$ and let $A'$ be the mutation of $A$ at $k$ (Definition~\ref{def:comp-mut}). Then $A'$ is an admissible quasi-Cartan companion of $\mu_k(\Gamma)=\Gamma'$ and $\Gamma'$ does not contain any subdiagram from Figure~\ref{fig:critical} as well.
\end{proposition}

We prove the proposition by obtaining a contradiction to the assumptions if any of the two stated properties is not true for $\mu_k(\Gamma)=\Gamma'$ as well. For this, first let us note that $A'$ is a quasi-Cartan companion of $\Gamma$ because $A$ is admissible.
Let $A''$ be the quasi-Cartan matrix obtained by mutating $A'$ at $k$. 
Then $A''$ is equal to $A$ up to a sign change at $k$, so $A''$ is an admissible quasi-Cartan companion of $\mu_k(\Gamma')=\Gamma$.  
We will obtain, in two lemmas, a contradiction to this or to the assumption that $\Gamma$ does not contain any diagram from Figure~\ref{fig:critical} if the conclusion of the proposition does not hold:

\begin{lemma}
\label{lem:mut class ext 1}
Let $\Gamma'$ be a diagram. Suppose that $A'$ is a quasi-Cartan companion of $\Gamma'$ which is semipositive of corank $1$ and let $A''$ be the quasi-Cartan matrix obtained by mutating $A'$ at $k$. Suppose also that $A'$ is \emph{not} admissible. Then either $A''$ is not an admissible quasi-Cartan companion of $\mu_k(\Gamma')=\Gamma$ or the diagram $\Gamma$ contains a subdiagram that belongs to Figure~\ref{fig:critical}. 
\end{lemma}

\noindent
{\bf Proof.} Since $A'$ is not admissible, there is a cycle $Z$ such that the restriction of $A$ on it is not admissible by Definition~\ref{def:admissible}. We first consider the case when $k$ is in $Z$. If $k$ is a source or sink of $Z$, then $\mu_k(Z)$ is also a cycle which does not satisfy the same condition of admissibility. If $k$ is not a source or sink (in $Z$), then either $A''$ is not a companion of $\mu_k(\Gamma')$ (this happens when $Z$ is a triangle) or the diagram obtained from $\mu_k(Z)$ by removing $k$ is a cycle which does not satisfy the same condition, so $A''$ is not admissible. 

We proceed by considering $k$ which is not in $Z$. By what we have just considered, we can assume that 

(*) the restriction of $A'$ to any cycle that contains $k$ is admissible. 

\noindent
For convenience, we will denote the subdiagram $\{Z,k\}$ by $Zk$. 
Note also that, since $A'$ is semipositive, the weight of any edge is at most $4$ (Proposition~\ref{prop:semidefinite comp}(i)).

{Case 1.} \emph{$Z$ is an oriented cycle.} 
Note that in this case $\prod (-A'_{i,j})$ over all edges of $Z$ is positive, so $Z$ is simply-laced by Proposition~\ref{prop:semidefinite comp}(iv). Then the restriction of $A$ to $Z$ has a non-zero radical vector $u$, which is a radical vector for $A'$ as well (Proposition~\ref{prop:semi-sub}). 
Applying some sign changes if necessary, we can assume that the coordinates of $u$ are equal to $1$ in the vertices of $Z$. Since $A'$ has corank $1$, the restriction of $A'$ to any subdiagram which does not contain $Z$ is positive. This implies, in particular, that any cycle $C$ which contains $k$ is oriented because if $C$ is non-oriented then, by Proposition~\ref{prop:double in admissible} (ii), the restriction of $A'$ to $C$ is not positive (this restriction is admissible by the assumption (*)). Similarly the weight of any edge which is adjacent to $k$ is less than $4$ (Proposition~\ref{prop:semidefinite comp}(ii)). 
If $k$ is connected to exactly one vertex in $Z$, then obviously $Z$ will be a subdiagram of $\Gamma$ such that the restriction of $A''$ to it is not admissible. Thus we can assume that $k$ is connected to at least two vertices in $Z$. 

Let us assume that $k$ is connected to $Z$ by an edge whose weight is $w=1,2,3$. Then, by the definition of a diagram, any edge connecting $k$ to $Z$ has weight $w$ respectively. We note that if $k$ is connected to two vertices in $Z$ which are not adjacent, then there is a non-oriented cycle that contains $k$ (because $Z$ is oriented), contradicting our assumptions. 
Thus we can assume that $k$ is connected to exactly two vertices, say $z_1,z_2$, in $Z$ and $z_1,z_2$ are adjacent; then note that the restriction of $A'$ to the edges $\{k,z_1\}$ and $\{k,z_2\}$ have opposite signs (so that $u$ is a radical vector). If $w=2,3$, then the effect of $\mu_k$ on $Z$ is to reverse the edge $\{z_1,z_2\}$ so that in $\mu_k(\Gamma')$, the subdiagram on $Z$ is a non-oriented cycle and the restriction of $A''$ to it is not admissible, contradiction.
(In fact, here, it is enough to take $w=2$ because if $w=3$, then the restriction of $A'$ on the subdiagram $\{z_1,z_2,k\}$ is not positive, contradicting our assumptions.) If $w=1$, then the effect of $\mu_k$ on $Z$ is to destroy the edge $\{z_1,z_2\}$ so that in $\mu_k(\Gamma')$, the subdiagram $\mu_k(Zk)$ is an oriented cycle and the restriction of $A''$ to it is not admissible. 

{Case 2.} \emph{$Z$ is a non-oriented cycle.} 
Note that in this case $\prod (-A'_{i,j})$ over all edges of $Z$ is negative. 
Let us first assume that $k$ is connected to a vertex $z$ in $Z$ by an edge $e$ whose weight is $4$. Let $z_1,z_2$ be the vertices which are adjacent to $z$ in $Z$. Then, by Proposition~\ref{prop:semidefinite comp}(iii) and Proposition~\ref{prop:double in admissible} (ii), the vertex $k$ is connected to both $z_1,z_2$ such that the triangles $T_1=\{k,z,z_1\}$ and $T_2=\{k,z,z_2\}$ are oriented, and $k$ is not connected to any other vertex on $Z$. Then both edges $\{k,z_1\}$ and $\{k,z_2\}$ have the same orientation, thus there is a non-oriented cycle $C$ which contains the edges  $\{k,z_1\}$ and $\{k,z_2\}$ (together with the edges from $Z$ which are not adjacent to $z$). By our assumption (*), the restriction of $A'$ to $T_1$ and $T_2$ is admissible; this implies that $\prod (-A'_{i,j})$ over all edges of $C$ is negative, so the restriction of $A'$ to $C$ is not admissible, contradicting (*). 

We now consider subcases assuming that the weight of any edge connecting $k$ to $Z$ is less than $4$. 


{Subcase 2.1.} \emph{$k$ is connected to exactly one vertex in $Z$.} 
Then obviously $Z$ will be a subdiagram of $\Gamma$ such that the restriction of $A''$ to $Z$ is not admissible.

{Subcase 2.2.} \emph{$k$ is connected to exactly two vertices in $Z$.} 
Say $k$ is connected to $z_1$ and $z_2$. Let us first assume that $Z$ contains an edge whose weight is equal to $4$. Then, by Proposition~\ref{prop:semidefinite comp}(iii), $Z$ is a triangle such that the weight of edge $e=\{z_1,z_2\}$ is equal to $4$. Furthermore, the edges $\{k,z_1\}$ and $\{k,z_2\}$ have equal weights, say $w$, such that the triangle $\{k,z_1,z_2\}$ is oriented. Then $w=1$ or $w=2$ because if $w$ is equal to $3$ then $\Gamma'$ contains a subdiagram of type $G_2^{(1)}$, implying that $A'$ has corank greater than or equal to two. (Note that $w\ne 4$ by our assumption above). Similarly, if $w=2$ then the weights of the edges of $Z$ are $4,1,1$. Then we have the following: if $w=2$, then the effect of $\mu_k$ on $Z$ is to destroy the edge $e=\{z_1,z_2\}$ so that $\mu_k(Zk)$ is a non-oriented cycle such that the restriction of $A''$ to it is not admissible; if $w=1$, then in $\mu_k(\Gamma')$, $Z$ stays as a non-oriented cycle but the weight of the edge $e=\{z_1,z_2\}$ is replaced by $1$ keeping the sign of the corresponding entry of the companion, so the restriction of $A''$ to $Z$ is not admissible, thus $A''$ is not admissible. Thus for the rest of this subcase, we can assume that $Z$ does not contain any edge whose weight is equal to $4$.

{Subsubcase 2.2.1.} \emph{$z_1$ and $z_2$ are connected.} 
First let us assume that the triangle $T=\{k,z_1,z_2\}$ is non-oriented. By our assumption (*), the restriction of $A'$ to this triangle is admissible, so it is simply laced (Proposition~\ref{prop:double in admissible}(ii)). If $k$ is a source or sink of $T$, then by the definition of mutation, $Z$ will be a subdiagram of $\Gamma$ and the restriction of $A''$ to it is still not admissible. If $k$ is not a source or sink of $T$, then in $\mu_k(\Gamma')$, $Z$ stays as a non-oriented cycle but the weight of the edge $\{z_1,z_2\}$ is replaced by $4$ keeping the sign of the corresponding entry of the companion, so the restriction of $A''$ to $Z$ is not admissible, thus $A''$ is not admissible.
Let us now assume that the triangle $T=\{k,z_1,z_2\}$ is oriented. 
Then the effect of $\mu_k$ on $Z$ is either to destroy the edge $e=\{z_1,z_2\}$ or to reverse it. If $\mu_k$ destroys $e$, then in $\mu_k(\Gamma')$ the subdiagram $\mu_k(Zk)$ is a non-oriented cycle such that the restriction of $A''$ to it is not admissible. Let us now assume that $\mu_k$ reverses $e$. Then in $\mu_k(\Gamma')$ the subdiagram on $Z$ is a cycle and $\prod (-A''_{i,j})$ over all edges of $Z$ is positive, so $Z$ is a simply-laced non-oriented cycle in $\mu_k(\Gamma')$ (otherwise $A''$ is not admissible or not semipositive). In particular, the weights of the edges $\{k,z_1\}$ and $\{k,z_2\}$ are equal. Also $Z$ has a 
vertex $v\ne z_1,z_2$ such that $v$ is a source or sink in $Z$ (because otherwise reversing $e$ produces an oriented cycle, contradicting that $Z$ is non-oriented in $\Gamma'$). Then we have the following: if $k$ is connected to $z_1$ (and $z_2$) by an edge of weight $2$, then $\mu_k(Zk)$ is of type $\check{B} ^{(1)}_{n}(r)$; 
if $k$ is connected to $Z$ by an edge of weight $3$, then in $\mu_k(Zk)$ the edges $\{k,z_1\}$ and $\{k,z_2\}$ are contained in seperate subdiagrams of type $G_{2}^{(1)}$, this implies that $A''$ has corank at least two (Proposition~\ref{prop:semi-sub}), contradicting our assumption.

{Subsubcase 2.2.2.} \emph{$z_1$ and $z_2$ are not connected.} In $Zk$ there are two cycles, say $C_1,C_2$, that contain $k$. By Lemma~\ref{lem:C positive} and (*), one of these cycles, say $C_1$, is non-oriented, so it is simply-laced (Proposition~\ref{prop:semidefinite comp}(iv)). Thus any edge connecting $k$ to $Z$ has weight $1$. Also by  Proposition~\ref{prop:semidefinite comp}(vi), the cycle $C_2$ is an oriented square. Given all this, let us note that the cycle $C_1$ has a source or sink which is not connected to $k$ because otherwise $Z$ needs to be oriented. Now we have the following: if $C_2$ is simply-laced, then $\mu_k(Zk)$ is of type $\check{D} ^{(1)}_{n}(r)$; if $C_2$ is not simply-laced, then it contains a subdiagram $S$ of type $C_{2}^{(1)}$ or $G_{2}^{(1)}$ such that $k$ is not in $S$, so there is a sincere radical vector for the restriction of $A'$ to $S$, which is also a radical vector for $A'$ (Proposition~\ref{prop:semi-sub}). Then $A'$ has corank $\geq 2$ because the restriction of $A'$ to $C_1$ also has a sincere radical vector. This contradicts the assumption of the proposition.


{Subcase 2.3.} \emph{$k$ is connected to exactly three vertices in $Z$.} If $Z$ contains an edge whose weight is equal to $4$, then the subcase is treated by similar arguments as in the Subcase 2.2. above. Let us assume that $Z$ does not contain any edge whose weight is equal to $4$.
In $Zk$ there are three cycles, say $C_1,C_2,C_3$, that contain $k$. One of these cycles, say $C_1$, is non-oriented, so simply-laced (by Proposition~\ref{prop:double in admissible}(ii) and (*), note that the restriction of $A'$ to $C_1$ is not positive). If $C_2$ or $C_3$ has more than $3$ vertices, then it contains a vertex which is connected to exactly one vertex in $C_1$, contradicting semipositiveness of $A'$ (Proposition~\ref{prop:semidefinite comp}(vi)). Thus we can assume that $C_2$ and $C_3$ are triangles. Let us denote by $v$ the vertex in $Z$ which is common to $C_2$ and $C_3$.
Now we have the following: if $C_2,C_3$ (so $Zk$) are simply-laced then $v$ is connected to exactly an odd number of vertices in $C_1$, contradiction (Proposition~\ref{prop:semidefinite comp}(vii)); if $C_2,C_3$ are not simply-laced, then they are oriented (Proposition~\ref{prop:double in admissible}(ii)) and the weights of the edges connecting $v$ to $C_1$ are equal (by the definition of a diagram), so, in $\mu_k(Zk)$, the vertex $v$ is connected to exactly one vertex in the non-oriented cycle $\mu_k(C_1)$ (note that $k$ is a source or a sink in $C_1$ so $\mu_k(C_1)$ is also a non-oriented cycle), contradiction by Proposition~\ref{prop:semidefinite comp}(vi).

{Subcase 2.4.} \emph{$k$ is connected to exactly four vertices in $Z$.} 
In this subcase there are four cycles, say $C_1,C_2,C_3,C_4$, that contain $k$. One of these cycles say $C_1$ is non-oriented, so simply-laced (by Lemma~\ref{lem:C positive} and (*); note that the restriction of $A'$ to $C_1$ is not positive). Then the restriction of $A'$ to each of $C_2,C_3,C_4$ is positive (otherwise $A'$ has higher corank by Proposition~\ref{prop:semi-sub}), so they are oriented (note then that $k$ is not a source or sink in $C_1$). Suppose that $C_2,C_3$ are adjacent to $C_1$. If any of $C_2,C_3$ has more than $3$ vertices, then it contains a vertex connected to exactly one vertex in $C_1$, contradicting semipositiveness of $A'$ by Proposition~\ref{prop:semidefinite comp}(vi). Thus we can assume that $C_2,C_3$ are (oriented) triangles.

Under all these assumptions, if the subdiagram $Zk$ is simply-laced, then we have the following: if $C_1$ has more than three and $C_4$ has three vertices, then $\mu_k(Zk)$ is of type $\check{D}^{(1)}_{n}(1,r)$; if each $C_1$ and $C_4$ has more than three vertices, then $\mu_k(Zk)$ is of type  $\check{D}_{n}^{(1)}(1,r,s)$; if each $C_1$ and $C_4$ has exactly three vertices, then $\mu_k(Zk)$ is of type $\check{D}^{(4)}_{n}$; if $C_1$ has exactly three vertices and $C_4$ has more, then $\mu_k(Zk)$ is of type $\check{D}^{(4)}_{n}(1,r)$. 

Let us now assume that $Zk$ is not simply-laced. If $Zk$ has an edge whose weight is equal to $3$, then it contains a subdiagram of type $G_2^{(1)}$, implying that $A'$ has corank greater than or equal to two (Proposition~\ref{prop:semi-sub}). For the same reason, $Zk$ does not contain any edge whose weight is equal to $4$. Thus the weight of any edge is $1$ or $2$. Let us note that, by the definition of a diagram, if any of the oriented triangles $C_2$ or $C_3$ is not simply-laced, then all $C_2,C_3,C_4$ are not simply-laced. Thus, in any case, the cycle $C_4$ is not simply-laced, therefore it is oriented by Proposition~\ref{prop:double in admissible}(ii) (note that the restriction of $A'$ to $C_4$ is admissible by our assumption (*)). This implies that the vertex $k$ is neither a source nor a sink of the non-oriented (simply-laced) cycle. Therefore if $C_2$ or $C_3$ is not simply-laced, then $\mu_k(Zk)$ contains a subdiagram of type $\check{B}^{(4)}_{3}$ or $\check{B}_{l}^{(1)}(1,r)$, for some $l \geq 4$; if $C_2$ and $C_3$ are simply-laced, then $Z$ (and $C_4$) contains a subdiagram of type $C_{l}^{(1)}$ for some $l$.  This implies that $A'$ has corank at least $2$ by Proposition~\ref{prop:double in admissible}(iv) and Proposition~\ref{prop:semi-sub}, which is a contradiction.

{Subcase 2.5.} \emph{$k$ is connected to at least five vertices in $Z$.} Then, in $Zk$, there are at least five cycles that contain $k$. Let 
us first assume that $k$ is connected to an odd number of vertices in $Z$. Then there is a non-oriented cycle $C\subset Zk$ which contains $k$. By
Proposition~\ref{prop:double in admissible}(ii) and (*), the cycle $C$ is simply-laced. There is a vertex in $Z$ which is connected to exactly one
vertex (which is $k$) in $C$. Then, by Proposition~\ref{prop:semidefinite comp}(vi), the companion $A'$ is indefinite, contradicting the assumption of the lemma. 

Let us now assume that $k$ is connected to an even number of vertices in $Z$. By Lemma~\ref{lem:C positive}, $k$ is contained in a cycle $C\subset Zk$ such that the product $\prod(-A_{i,j})$ over all edges of $C$ is positive. This implies that $C$ is non-oriented because the restriction of $A'$ to $C$ is admissible by our assumption (*).  Also there is a vertex in $Z$ which is connected to exactly one vertex (which is $k$) in $C$. If $C$ is simply-laced, then the companion $A'$ is indefinite by Proposition~\ref{prop:semidefinite comp}(vi), contradicting the assumption of the lemma. If $C$ is not simply-laced, the same contradiction is provided by Proposition~\ref{prop:double in admissible}(ii).

The proof of Lemma~\ref{lem:mut class ext 1} is completed.

To proceed with the proof of Proposition~\ref{prop:mut class ext}, we can now assume that $A'$ is admissible. To complete the proof, we need to show that $\Gamma'$ does not contain any diagram from Figure~\ref{fig:critical}. We show this by obtaining a contradiction:

\begin{lemma}
\label{lem:mut class ext 2}
Suppose that $\Gamma '$ is a diagram and let $A'$ be an admissible quasi-Cartan companion which is semipositive of corank $1$. Suppose also that 
$\Gamma '$ contains a subdiagram $X$ that belongs to Figure~\ref{fig:critical}. Let $A''$ be the mutation of $A'$ at a vertex $k$. Then either $A''$ is not an admissible quasi-Cartan companion of $\mu_k(\Gamma')=\Gamma$ or $\Gamma$ contains a subdiagram that belongs to Figure~\ref{fig:critical}. 
\end{lemma}
\noindent
{\bf Proof.}
If $k$ is on $X$, then the lemma follows from a direct check. Then, to consider $k$ which is not on $X$, we can assume, by  Proposition~\ref{prop:double in admissible}(v), that 

(**) $k$ is not contained in any subdiagram from Figures~\ref{fig:extended-dynkin-diagrams},~\ref{fig:minimal} and \ref{fig:critical}

\noindent
because $X$ already contains an edge of weight $4$ or a non-oriented cycle. In particular, we assume that any cycle that contains $k$ is oriented. 
(In fact, we can assume that $k$ is not contained in any subdiagram $M$ of minimal infinite type, because any admissible companion of $M$ is    semipositive of corank $1$ with a sincere radical vector, see Theorem~\ref{th:minimal}). 
For convenience, we denote the subdiagram $\{X,k\}$ by $Xk$. If $k$ is connected to exactly one vertex in $X$, then $X$ is also a subdiagram of $\Gamma$. Thus, for the rest of the proof, we can assume that $k$ is connected to at least two vertices in $X$. We assume that the vertices of $X$ are labeled as in Figure~\ref{fig:critical}

{Case 1.} \emph{$X$ is of type $\check{D}^{(1)}_{n}(r)$.} Let us first assume that $k$ is connected to $X$ by an edge of weight $2$ or $3$. Then, by the definition of a diagram, any edge connecting $k$ to $X$ has weight $2$ or $3$ respectively. 
Thus we have the following: if $k$ is connected to two vertices $x_1$ and $x_2$ which are not connected in $X$, then 
the subdiagram $\{k,x_1,x_2\}$ is of type ${C}_2^{(1)}$; otherwise, it can be checked easily that $k$ is contained in a subdiagram of type ${B_l}^{(1)}$ for some $l$ or ${G}^{(1)}_{2}$, contradicting (**). 

We proceed by considering the cases where any edge connecting $k$ to $X$ has weight $1$ (so $Xk$ is simply-laced).

{Subcase 1.1.} \emph{$k$ is connected to both $b_1$ and $b_2$.} Then $k$ is not connected to any of $a_1$ or $c_1$, 
because otherwise there would be a non-oriented triangle that contains $k$. 
Then the subdiagram $\{b_1,b_2,a_1,c_1,k\}$ is of minimal infinite type $\check{D}_4^{(1)}(3)$, contradicting (**).

{Subcase 1.2.} \emph{$k$ is connected to only one of $b_1,b_2$.} Say $k$ is connected to $b_2$. 
Note that $k$ is also connected to another vertex among  $b_3,...,b_r$ (Proposition~\ref{prop:semidefinite comp}(vii)).

Let us first consider the subcase where $k$ is not connected to any of $a_1$ and $c_1$. 
If $k$ is not connected to $b_3$, then the subdiagram $\{a_1,c_1,b_2,b_3,k\}$ is of type ${D}^{(1)}_{4}$, contradicting (**). Let us now assume that $k$ is connected to $b_3$. If $k$ is not connected to any other $b_i$, then $\mu_k(Xk)$ is of type $\check{D}^{(1)}_{n+1}(r+1)$  (note that here the subdiagram  $\{k,b_2,b_3\}$ is oriented by (**)) . If $k$ is connected to $b_i$ such that $i>3$, then we can  assume, without loss of generality, that $k$ is not connected to any  $b_j$ for $j>i$. Then either 
the subdiagram $\{k,b_i,b_{i+1},...,b_r,b_1,b_2\}$ is a non-oriented cycle or 
the subdiagram $\{k,b_i,b_{i+1},...,b_r,b_1,b_2,a_1,c_1\}$ is of type ${D}^{(1)}(r-i+4)$, contradicting (**). 

Let us now consider the subcase where $k$ is connected to $a_1$ or $c_1$.
If $k$ is connected to both of them, then the cycle $\{k,a_1,b_1,c_1\}$ is non-oriented, so assume without loss of generality that $k$ is connected only to $a_1$. Let $b_i$, $i\geq 3$, be the vertex such that $k$ is connected to $b_i$ but not connected to any $b_j, j>i$. Then the subdiagram $\{k,b_2,a_1\}$ or $\{k,b_i,b_{i+1},...,b_r,b_1,a_1\}$ is a non-oriented cycle, contradicting (**).

{Subcase 1.3.} \emph{$k$ is not connected to any of $b_1,b_2$.}
Let us first assume that $k$ is connected to $a_1$ or $c_1$, say connected to $a_1$. If $k$ is connected to $c_1$ as well, then the cycle $\{k,a_1,c_1,b_1\}$ is non-oriented. If $k$ is not connected to $c_1$, then it is connected to a vertex $b_i$, $3 \leq i\leq r$, and so
there are two cycles $C_1,C_2$ that contain the edge $\{k,a_1\}$ together with one of the edges $\{a_1,b_1\}$ or $\{a_1,b_2\}$. One of the cycles $C_1,C_2$ is non-oriented because the triangle $\{a_1,b_1,b_2\}$ is oriented, contradicting (**). 
If $k$ is not connected to any of $a_1$ or $c_1$, then, by the same argument in Subcase 1.2 above, either $\mu_k(Xk)$ is of type $\check{D}^{(1)}_{n+1}(r+1)$ or $k$ is contained in a subdiagram which is a non-oriented cycle or is of type ${D}^{(1)}(r-t)$ for some $t<r$, contradicting (**).  

{Case 2.} \emph{$X$ is of type $\check{D}^{(1)}_{n}(m,r)$.} 
As in Case 1 above, if $k$ is connected to $X$ by an edge of weight $2$ or $3$, then $k$ is contained in a subdiagram of type ${C}^{(1)},{B}^{(1)}$ or ${G}^{(1)}_{2}$, contradicting (**). 
Thus, for the rest of this case, we assume that any edge connecting $k$ to $X$ has weight $1$. We denote the non-oriented cycle in $X$ by $C$.

{Subcase 2.1.} \emph{$k$ is connected to $C$.} By Proposition~\ref{prop:semidefinite comp}(vii), the vertex $k$ is connected to an even number of vertices in $C$. 
Let us first assume that 
$k$ is not connected to any $a_i$, $i=1,...,m$ nor to $c_1,c_2$. Let $C_1,...,C_r$ be the (oriented) cycles that contain $k$. If one of these cycles, say $C_i,$ contains the edge $\{b_1,b_2\}$, then the subdiagram $\{C_i,a_1,...,a_m,c_1,c_2\}$ is of type ${D}^{(1)}(m,t)$ for some $t\leq r$, contradicting (**). If such a cycle does not exist, then $k$ is connected to exactly two vertices, say $b_i,b_j$, in $C$ which are connected and $\{b_i,b_j\} \ne \{b_1,b_2\}$. 
Then $\mu_k(Xk)$ is of type $\check{D}^{(1)}(m,r+1)$. Let us now assume that $k$ is connected to $a_j$ or $c_1,c_2$; we can assume without loss of generality that $k$ is not connected to $a_{i}$, $i<j$ (take $j=m+1$ if $k$ is not connected to any $a_i$). 
Then, since $k$ is connected to an even number of vertices in $C$, there are two cycles $C_1,C_2$ that contain the edge $\{k,a_j\}$ together with one of the edges $\{a_1,b_2\}$ or $\{a_1,b_1\}$. Since the triangle $\{a_1,b_1,b_2\}$ is oriented, one of the cycles $C_1,C_2$ is non-oriented, contradicting (**).

{Subcase 2.2.} \emph{$k$ is not connected to $C$.}
Let us first note that if $k$ is not connected to any of $a_1,...,a_m$, then it is connected to both $c_1,c_2$, so $\mu_k(Xk)$ is of type $\check{D}^{(1)}_{n+1}(m,r,3)$. Let us now assume that $k$ is connected to $a_{i_1},....,a_{i_j}$, $1\leq j \leq m$, $1\leq {i_1}<....< {i_j}\leq m$. We note that if $i_2\ne i_1+1$, then the subdiagram $\{C,a_{1},....,a_{i_1},a_{i_1+1},k\}$ is of type $\check{D}^{(1)}(i_1,r)$; if  $i_2=i_1+1$ but $j \geq 3$, then the subdiagram $\{C,a_{1},....,a_{i_1},a_{i_2},...,a_{i_3},k\}$ is of type $\check{D}^{(1)}(i_1,r,i_3-i_2+2)$. Now there remain two subcases to consider. The first subcase is when $j=2$ such that $i_2=i_1+1$: if $k$ is connected to $c_1$ or $c_2$, say to $c_1$, then the subdiagram obtained from $Xk$ by removing $c_2$ is of type $\check{D}^{(1)}(i_1,r,m-i_2+3)$, contradicting (**); otherwise $\mu_k(Xk)$ is of type $\check{D}^{(1)}(m+1,r)$. 
Now the only subcase left is when $j=1$. If $i_1\ne m$, then the subdiagram $\{C,a_1,...,a_{i_1},a_{i_1+1},k\}$ is of type $\check{D}^{(1)}(i_1,r)$, contradicting (**). If $i_1=m$, then we have the following: if $k$ is not connected to one of $c_1,c_2$, say not connected to $c_2$, then the subdiagram obtained from $Xk$ by removing $c_1$ is of the same type as $X$, contradicting (**); if $k$ is connected to both, then $\mu_k(Xk)$ is of type $\check{D}^{(1)}(m+1,r)$.

{Case 3.} \emph{$X$ is of type $\check{D}^{(1)}_{n}(m,r,s)$.} Let us note that $X$ is very similar to the diagram $\check{D}^{(1)}_{n}(m,r)$, which we considered in Case 2 above. The case follows by similar arguments as in Case 2.



{Case 4. }\emph{$X$ is of type $\check{D}^{(4)}_{n}$.} 
We denote by $e$ the edge $\{b_1,b_2\}$ whose weight is $4$. 

{Subcase 4.1.} \emph{$k$ is connected to $e$.} Note that the subdiagram $\{e,k\}$ is an oriented triangle (by Proposition~\ref{prop:double in admissible}(i)). If $k$ is connected to a vertex which is not adjacent to $e$, then by the same argument as in Subcase 2.1, there is a non-oriented cycle that contains $k$, contradicting (**). Therefore we can assume that $k$ is not connected to any vertex other than $b_1$ and $b_2$. If $k$ is connected $b_1$ and $b_2$ by an edge of weight $2$ or $3$, then $k$ is contained in a subdiagram of type ${B}^{(1)}$ or ${G}^{(1)}_{2}$ respectively, contradicting (**), otherwise $\mu_k(Xk)$ is of type $\check{D}^{(1)}(m,r)$, with $r=3$, $m=n-3$.


{Subcase 4.2.} \emph{$k$ is not connected to $e$.} 
The subcase follows by similar arguments as in Subcase 2.2 above.

{Case 5.} \emph{$X$ is of type $\check{D}^{(4)}_{n}(m,r)$.} Let us note that $X$ is very similar to the diagrams in Cases 2 and 4. This case also  
follows by similar arguments as in these cases.


{Case 6.} \emph{$X$ is one of types $\check{B}^{(4)}_{n}$, $\check{B}^{(1)}_{n}(m,r)$ or $\check{B}^{(1)}_{n}(r)$.} Let us note that these diagrams are very similar to the diagrams $\check{D}^{(4)}_{n}$,  $\check{D}^{(1)}_{n}(m,r)$, $\check{D}^{(1)}_{n}(r)$ respectively. This case also follows by similar arguments as in these cases. 


Let us now prove the theorem:

\noindent
{\bf Proof of Theorem~\ref{th:mut class ext}} If $\Gamma(B)$ is an extended Dynkin diagram then it does not contain any diagram from  Figure~\ref{fig:critical} and $B$ has an admissible quasi-Cartan companion which is semipositive of corank $1$. 
The same conclusion holds for any skew-symmetrizable matrix whose diagram is mutation-equivalent to an extended Dynkin diagram by Proposition~\ref{prop:mut class ext}.





To prove the converse, let us assume that $B$ has an admissible quasi-Cartan companion $A$ which is semipositive of corank $1$ and $\Gamma=\Gamma(B)$ 
does not contain any diagram from Figure~\ref{fig:critical}. We will show that $\Gamma$ is mutation-equivalent to  an extended Dynkin diagram. Since $A$ is not positive, the diagram $\Gamma$ is not of finite type (Theorem~\ref{th:2-finite-class positive}), so it is mutation-equivalent to a diagram $\Gamma'$ which has an edge $e$ whose weight is $4$. Furthermore $\Gamma'$ has an admissible quasi-Cartan companion $A'$ which is semipositive of corank $1$ and it does not contain any diagram from Figure~\ref{fig:critical} (Proposition~\ref{prop:mut class ext}). Also, by Proposition~\ref{prop:double in admissible}(v), the diagram $\Gamma'$ does not contain any diagram from Figure~\ref{fig:extended-dynkin-diagrams} except the edge $e$ or from  Figure~\ref{fig:minimal} (in particular does not contain any non-oriented cycle).

We note that if a vertex $v$ is connected to $e$, then the subdiagram on $v,e$ is an oriented triangle (Proposition~\ref{prop:double in admissible}(i)). For any such $v$, we denote by $Pv$ the subdiagram on vertices which are connected to $v$ by a path that does not contain any vertex which is adjacent to $e$. Let us  denote the vertices connected to $e$ by $v_1,v_2,...,v_r$. For any $v_i\ne v_j$ connected to $e$, the subdiagrams $Pv_i$ and $Pv_j$ are disjoint because otherwise there is a non-oriented cycle in $\Gamma'$, contradicting our assumption. Thus any path connecting a vertex in $Pv_i$ to $Pv_j$, $i\ne j$, contains a vertex which is adjacent to $e$. 


Let us first consider the case where $\Gamma'$, so $\Gamma$, represents a skew-symmetric matrix. Recall that $\Gamma'$ does not contain any diagram from Figure~\ref{fig:critical}, in particular it does not contain any subdiagram of type $\check{D}^{(4)}_{n}$ or $\check{D}^{(4)}_{n}(m,r)$, therefore for any $v$ connected to $e$ the subdiagram $Pv$ does not contain any subdiagram which is of type $D_4$ or formed by two adjacent cycles. This implies that  $Pv$ is mutation-equivalent to $A_n$ \cite[Corollary 5.15]{S2}; applying some mutations if necessary, we can assume that $Pv$ is of type $A_n$ such that $v$ is the end vertex of $Pv$ (otherwise there is a subdiagram of type $\check{D}^{(4)}_{4}$; also note that if mutations are applied, then the resulting diagram also has an admissible quasi-Cartan companion which is semipositive of corank $1$ and it does not contain any diagram from Figure~\ref{fig:critical} by Proposition~\ref{prop:mut class ext}, so we will not lose any generality). 
Then $r\leq 3$ because otherwise there is a subdiagram of type $D_4^{(1)}$, which belongs to Figure~\ref{fig:extended-dynkin-diagrams}, contradicting our assumption. If $r\leq 2$, then $\Gamma'$ is mutation-equivalent to $A_n^{(1)}$, $n\geq 1$, as can be seen easily by applying mutations at the vertices which are connected to $e$. Let us now assume that $r=3$. If all $Pv_1,Pv_2,Pv_3$ have at least two vertices, then there is a subdiagram of type $E_6^{(1)}$, which contradicts our assumption, so we can assume that $Pv_1$ has exactly one vertex (which is $v_1$). Similarly, if $Pv_2$ and $Pv_3$ both have at least three vertices, then there is a subdiagram of type $E_7^{(1)}$, so we can assume that $Pv_2$ has at most two vertices. If $Pv_2$ has exactly one vertex (which is $v_2$), then $\Gamma'$ is mutation-equivalent to $D_n^{(1)}$. If $Pv_2$ has exactly two vertices, then $Pv_3$ has at most four vertices (otherwise there is a subdiagram of type $E_8^{(1)}$), then we have the following: if $Pv_3$ has exactly one vertex (which is $v_3$), then $\Gamma'$ is mutation-equivalent to $D_5^{(1)}$; if it has exactly two vertices, then $\Gamma'$ is mutation-equivalent to $E_6^{(1)}$; if it has exactly three vertices, then $\Gamma'$ is mutation-equivalent to $E_7^{(1)}$; if it has exactly four vertices, then $\Gamma'$ is mutation-equivalent to $E_8^{(1)}$.


Let us now consider the case where $\Gamma'$ does not represent a skew-symmetric matrix, so $\Gamma'$ has an edge whose weight is $2$ or $3$. Then any such edge of weight $2$ or $3$ is connected to $e$, because otherwise $\Gamma'$ contains a subdiagram of the following types $C_n^{(1)}$, $\check{B}^{(4)}_{n}$, $G_2^{(1)}$ or a three-vertex tree $T$ with edge-weights $2,3$, contradicting our assumption (the first three types belong to Figure~\ref{fig:extended-dynkin-diagrams} or Figure~\ref{fig:critical}; the restriction of $A'$ to $T$ is indefinite). For the same reasons, there is exactly one vertex, say $v_1$, which is connected to $e$ by an edge of weight $2$ or $3$. Then, note in particular, that the only edges in $\Gamma'$ whose weights are $2$ or $3$ are the two edges that connect $v_1$ to $e$. If $v_1$ is connected to $e$ by an edge of weight $3$, then $\Gamma'$ does not contain any other vertex (because otherwise there is a subdiagram of type $G_2^{(1)}$), so $\Gamma'$ is mutation-equivalent to $G_2^{(1)}$. 

Thus for the rest of the proof we can assume that the weight of any edge connecting $v_1$ to $e$ is $2$. As in the skew-symmetric case above, for any $v_i$ connected to $e$, the subdiagram $Pv_i$ does not contain any subdiagram which is of type $D_4$ or $B_n^{(1)}$ or formed by two adjacent cycles. This implies that each $Pv_i$ is mutation-equivalent to $A_n$ \cite[Corollary 5.15]{S2}; applying some mutations if necessary, we can assume that $Pv_i$ is of type $A_n$ such that $v_i$ is the end vertex of $Pv_i$ (otherwise there is a subdiagram of type $\check{D}^{(4)}_{4}$ or $B_3^{(1)}$). 
Let us note that we have $r \leq 2$, because otherwise there is a subdiagram of type $B_3^{(1)}$. If $r=1$, then $\Gamma'$ is mutation-equivalent to $C_n^{(1)}$. If $r=2$, then $Pv_1$ has at most two vertices (because otherwise there is a subdiagram of type $F_4^{(1)}$) so we have the following: if $Pv_1$ has exactly one vertex, then $\Gamma'$ is mutation-equivalent to $B_n^{(1)}$; if $Pv_1$ has two vertices, then $\Gamma'$ is mutation-equivalent to $F_4^{(1)}$. This completes the proof of Theorem~\ref{th:mut class ext}.

\subsection{Proof of Theorem~\ref{th:ext skew}}
\label{sec:ext skew}

The implication (2) $\Rightarrow$ (1) trivially follows from the definition of a quasi-Cartan companion.  To show (1) $\Rightarrow$ (2), let us suppose that $X$ is in $\mathcal{S}$ such that $\Gamma=\mu_{i_r}...\mu_{i_1}(\Gamma(X))$ is an extended Dynkin diagram. Then $B=\mu_{i_r}...\mu_{i_1}(X)$ is a skew-symmetrizable matrix whose diagram is $\Gamma$. Then it follows from a direct check on Tables of \cite[Chapter 4]{K} that $B$ has a quasi-Cartan companion $A$ which is a generalized Cartan matrix (of affine type). (Note that $X,B$ and $A$ share the same (skew-)symmetrizing matrix $D$). To prove the uniqueness of $A$,  let us assume that $X$ is mutation-equivalent to a skew-symmetrizable matrix $B'$, say $B'=\mu_{j_s}...\mu_{j_1}(X)$, which has another generalized Cartan matrix $A'$ as a quasi-Cartan companion. Then $B'=\mu_{j_s}...\mu_{j_1}\mu_{i_1}...\mu_{i_r}(B)$. On the other hand, since $A$ and $A'$ are admisssible, by Proposition~\ref{prop:mut class ext}, $A'$ can be obtained from $A$ by the same sequence of mutations possibly with simultaneous sign changes in rows and columns. This implies, in particular, that $A$ and $A'$ are equivalent. Thus $\mathcal{S}$ determines $A$ uniquely.

For the converse, let $A$ be an affine type generalized Cartan matrix which is not of type $A_n^{(1)}$, $n\geq 2$. Let $B$ be any skew-symmetrizable matrix which has $A$ as a quasi-Cartan companion. Then note that for any such choice of $B$ its diagram is a tree diagram (so $A$ is an admissible quasi-Cartan companion). Also any two orientations of a tree diagram can be obtained from each other by a sequence of mutations (at source or sink vertices, i.e. by reflections), which implies that any two choices for $B$ are mutation-equivalent \cite[Proposition~9.2]{CAII}. Also, by our argument above via Proposition~\ref{prop:mut class ext}, another skew-symmetrizable $B'$ defined in the same way by a different affine type $A'$ is not mutation-equivalent to $B$ (otherwise $A$ and $A'$ are equivalent). Thus the mutation class of $B$ is uniquely determined by $A$.

Different non-cyclic orientations of a cycle are not necessarily mutation-equivalent to each other. For this reason, there are non-oriented cycles which are not mutation-equivalent while they have the same generalized Cartan matrix $A_n^{(1)}$, $n\geq 2$, as an admissible quasi-cartan companion. We refer to \cite{B} for a study of mutation classes of those diagrams.  


\subsection{Proof of Theorem~\ref{th:mut class ext char}}
\label{sec:mut class ext char}

By Proposition~\ref{prop:mut class ext}, any diagram which is mutation-equivalent to an extended Dynkin diagram has an admissible quasi-Cartan companion which is semipositive of corank $1$. For the converse, first it can be checked easily that any diagram from Figure~\ref{fig:critical} which corresponds to a skew symmetric matrix, (i.e. a diagram of type $\check{D}$) is mutation-equivalent to a diagram which contains $\check{D}^{(4)}_{4}$ (with $5$ vertices). Let us assume without loss of generality that $\check{D}^{(4)}_{4}$ is oriented in such a way that there are two edges oriented away from the vertex in the "center" and two edges oriented towards it. Then mutating at the "center" 
results in a diagram which does not have any admissible quasi-Cartan companion. Thus, if $\Gamma$ is the diagram of a skew symmetric matrix such that any diagram in its mutation class has an admissible companion which is semipositive of corank $1$, then it does not contain any subdiagram which belongs to Figure~\ref{fig:critical}, implying that $\Gamma$ is mutation-equivalent to an extended Dynkin diagram by Theorem~\ref{th:mut class ext}.



\subsection{Proof of Theorem~\ref{th:minimal}}
\label{sec:minimal}

To prove the first statement, suppose that $\Gamma$ has an admissible quasi-Cartan companion $A$ which is semipositive of corank $1$ with a sincere radical vector, so any non-zero radical vector is also sincere. Let $k$ be an arbitrary vertex of $\Gamma$. Let $\Delta$ be the subdiagram obtained from $\Gamma$ by removing $k$. Then the restriction $A'$ of $A$ to $\Delta$ is positive: otherwise $A'$ has a non-zero radical vector $u$, which is a radical vector for $A$ as well (see Proposition~\ref{prop:semi-sub}), however $u$ is not sincere, contradicting that any radical vector for $A$ is sincere. Thus $\Delta$ is of finite type (Theorem~\ref{th:2-finite-class positive}). Since $k$ is an arbitrary vertex, any subdiagram of $\Gamma$ is of finite type, so $\Gamma$ is of minimal infinite type (here $\Gamma$ is of infinite type because $A$ is not positive).



For the converse, let us recall that minimal infinite type diagrams have been computed explicitly in \cite{S2}: it follows from a direct check that each of them has an admissible quasi-Cartan companion which semipositive of corank $1$ with a sincere radical vector. (Applying sign changes if necessary, the coordinates of this radical vector can be assumed to be positive). Here, for a minimal infinite type diagram $\Gamma$ which corresponds to a skew-symmetric matrix, we offer an alternative proof: The statement is true for any simply-laced non-oriented cycle (Proposition~\ref{prop:semidefinite comp}(v)). Thus we can assume that $\Gamma$ does not have any non-oriented cycles. Then $\Gamma$ has an admissible quasi-Cartan companion $A$ \cite[Corollary~5.2]{BGZ}. Since any proper subdiagram of $\Gamma$ is of finite type, the restriction of $A$ to any proper subdiagram is positive (Theorem~\ref{th:2-finite-class positive}). This implies, by \cite[Theorem~2, Section~1.0]{R}, that the companion $A$ is semipositive of corank $1$ with a sincere radical vector (recall $\Gamma$ has at least three vertices). This completes the proof of the first statement. 

To prove the second part, let $\Gamma$ be a diagram of minimal infinite type. By the first part, it has an admissible quasi-Cartan companion which is semipositive of corank $1$ with a sincere radical vector. This implies that any admissible companion of $\Gamma$ has a sincere radical vector. 
On the other hand, any admissible companion of a diagram that belongs to Figure~\ref{fig:critical} has a non-zero radical vector which is not sincere
(Proposition~\ref{prop:semidefinite comp}(ii,v)). Therefore $\Gamma$ does not contain any diagram which belongs to Figure~\ref{fig:critical}.  Thus $\Gamma$ is mutation-equivalent to an extended Dynkin diagram by Theorem~\ref{th:mut class ext}. This completes the proof of the theorem. 


\noindent
{\bf Remark.}  Most of the minimal infinite type diagrams correspond to skew-symmetric matrices (i.e. most of them are quivers) and their quasi-Cartan companions as described in the theorem can be found in \cite{HV}. More explicitly, \cite{HV} gives a list of symmetric matrices (viewed as sign assignments on underlying graphs of quivers) that represent a class of quadratic forms which are called "Tits forms of tame concealed algebras"; those symmetric matrices turn out to be quasi-Cartan companions of minimal infinite type diagrams. The relation between minimal infinite type diagrams and tame concealed algebras in the setup of cluster categories have been studied in \cite{BRS}. 
\subsection{Proof of Theorem~\ref{th:acyclic fmc}}
\label{sec:acyclic fmc}

We prove the theorem using the following two lemmas, which give some basic types of diagrams whose mutation classes are  infinite.
\begin{lemma}\label{lem:infinite}
Let $\Gamma$ be a connected diagram which has at least three vertices.
\begin{enumerate}
\item[(i)] If $\Gamma$ has an edge whose weight is greater than $4$, then it has an infinite mutation class. 

\item[(ii)] Suppose that $\Gamma$ has exactly three vertices and has an edge whose weight is $4$. Then  $\Gamma$ has a  finite mutation class if and only if it is an oriented triangle with edge weights $4,1,1$ or $4,4,4$ or $4,2,2$ or $4,3,3$.

\item[(iii)] If $\Gamma$ is a non-simply-laced cycle which is non-oriented, then it has an infinite mutation class.

\item[(iv)] Suppose that $\Gamma$ does not have any edge whose weight is greater than or equal to $4$.  
If $\Gamma$ has a non-oriented cycle $C$ such that there is a vertex $k$ which is connected to exactly an odd number of vertices in $C$, then it has an infinite mutation class.

\item[(v)] Suppose that $\Gamma$ does not contain any oriented cycle but has at least two non-oriented cycles. Then $\Gamma$ has an infinite mutation class. 
\end{enumerate}
\end{lemma}

\noindent
Statements (i),(ii),(iii) easily follow from the definitions. Let us prove (iv). By part (iii), we can assume that $C$ is simply-laced. First we consider the case where $k$ is connected to exactly one vertex, say $c$, in $C$. Let us assume first that $C$ is a triangle. Applying a mutation at a source or sink of $C$ if necessary, we can assume that $c$ is a source or sink; mutating at the vertex which is neither a source or sink, we obtain a diagram which contains a three-vertex tree which has an edge whose weight is $4$; then part (ii) applies. Let us now assume that $C$ has more than $3$ vertices. Then, applying a mutation at a source or sink of $C$ if necessary, we can assume that there is a vertex $c'$ in $C$, $c\ne c'$, which is neither a source nor a sink in $C$. Then in $\mu_{c'}(\Gamma)$, the subdiagram  $C'$ obtained from $C$ by removing $c'$ is a non-oriented cycle and $k$ is connected to exactly one vertex in $C'$. Then the statement (iv) follows by induction.


Let us now consider the case where $k$ is connected to exactly three vertices in $C$. Then there are three cycles, say $C_1,C_2,C_3$, that contain $k$; one of them, say $C_1$, is necessarily non-oriented. If $C_1$ is not simply-laced then part (iii) applies, so we can assume that $C_1$ is simply-laced. This implies that any edge connecting $k$ to $C$ has weight $1$. If one of the cycles $C_2$ or $C_3$ has more than three vertices, then there is a vertex in that cycle connected to exactly one vertex in $C_1$, which is the case we have considered above. Thus we can further assume that $C_2$ and $C_3$ are triangles. Given all this, we proceed as follows. If $C$ has exactly three vertices, then the statement follows from a direct check. If $C$ has more than three vertices, then one of the cycles $C_1,C_2,C_3$ also has more than three vertices; since $C_2$ and $C_3$ are triangles, the cycle $C_1$ must have at least four vertices. If any of $C_2$ or $C_3$ is non-oriented, then there is a vertex in $C_1$ which is connected to exactly one vertex in that cycle, which is the case we considered above. Then the only subcase left to consider is the case where both $C_2$ and $C_3$ are oriented. Then, in $\mu_{k}(\Gamma)$, the subdiagram $\{C,k\}$ consists of a non-oriented cycle $C'$ that contains $k$ and an additional vertex which is connected to exactly one vertex in $C'$, which is again the case we have considered. To consider the case where $k$ is connected to at least five vertices in $C$, we note that in this case there is a non-oriented cycle $C'$ which contains $k$ and there is a vertex in $C$ connected to exactly one vertex in $C'$, which is a case we have considered.


To prove part (v), we can assume that any cycle in $\Gamma$ is simply-laced by part(iii). Let us now suppose that $C$ is a cycle with minimal number of vertices in $\Gamma$. There is a vertex $k$ which is not in $C$ but connected to $C$. If $k$ is connected to $C$ by an edge $e$ of weight $4$, then there is a three-vertex tree that contains $e$, so part (ii) applies; if $k$ is connected to $C$ by an edge $e$ of weight $2$ or $3$, then $k$ is connected to exactly one vertex in $C$ (because we assumed that any cycle in $\Gamma$ is simply-laced), then part (iv) applies. Thus we can assume that
any edge connecting $k$ to $C$ has weight $1$.  Then we have the following. If $k$ is connected to an odd number of vertices in $C$, then part (iv) applies. If $k$ is connected to an even number of vertices and $C$ is a triangle or a square, then the statement follows from a direct check; if $C$ has at least five vertices, then there is a non-oriented cycle $C'$ containing $k$ such that another vertex $r$ is connected to exactly an odd number of vertices in $C'$, so if $r$ is connected to $C'$ by an edge of weight $4$ then part (ii) applies, otherwise part (iv) applies. This completes the proof of the lemma.

\begin{lemma}\label{lem:C indefinite}
Suppose that $\Gamma$ is a diagram with an indefinite admissible quasi-Cartan companion $A$. Suppose also that $\Gamma$ contains a subdiagram $X$ which is either an edge of weight $4$ or a non-oriented cycle. Let $u$ be a non-zero radical vector for the restriction of $A$ to $X$ (i.e. $u$ is in the span of the standard basis vectors which correspond to the vertices in $X$ and $x^TAu=0$ for all $x$ in the same span.). If $u$ is not a radical vector for $A$, then $\Gamma$ has an infinite mutation class. In particular, the conclusion holds if $A$ is non-degenerate.
\end{lemma}
\noindent
To prove the lemma, we can assume that the weight of any edge is at most $4$ (Lemma~\ref{lem:infinite}(i)). We first show the lemma for the case where $X$ is an edge whose weight is $4$. Since $u$ is not a radical vector for $A$, there is a three-vertex subdiagram $Y$ containing $X$ such that $u$ is not a radical vector for the restriction of $A$ to $Y$. Since $A$ is admissible, the subdiagram $Y$ is not an oriented triangle with weights $4,1,1$ or $4,4,4$ or $4,2,2$ or $4,3,3$ (otherwise $u$ becomes a radical vector for the restriction of $A$ to $Y$ as well), so $Y$, thus $\Gamma$, has an infinite mutation class by Lemma~\ref{lem:infinite}(ii). 

Let us now show the lemma for the case where $X$ is a non-oriented cycle. By  Lemma~\ref{lem:infinite}(iii), we can assume that $X$ is simply-laced. We can also assume, applying sign changes if necessary, that the restriction of $A$ to any edge of $X$ is $-1$. As before, there is an additional vertex $k$ which is connected to $X$ such that the restriction of $A$ to the subdiagram $Y=\{X,k\}$ does not have $u$ as a radical vector. We first consider the subcase where $k$ is connected to a vertex, say $z$, in $X$ by an edge $e$ whose weight is $4$. Let $z_1,z_2$ be the vertices which are connected to $z$ in $X$. Then $k$ is contained in a three-vertex subdiagram which is not as in Lemma~\ref{lem:infinite}(ii) unless the following holds: $k$ is connected to both $z_1$ and $z_2$ with edges of weight $1$ and $k$ is not connected to any other vertex in $X$ such that both triangles $\{k,z,z_1\}$ and $\{k,z,z_2\}$ are oriented; then, however, $u$ is a radical vector for the restriction of $A$ to the subdiagram $Y$, contradicting our assumption. Thus $\Gamma$ has an infinite mutation class. 
Let us now consider the remaining subcase, where all edges connecting $k$ to $X$ have weight less than $4$. Then such edges all have the same weight (because of the definition of a diagram), thus the number of those edges assigned $(-)$ is different from the ones assigned $(+)$ (not to have $u$ as a radical vector). Then, either $k$ is connected to an odd number of vertices in $X$, so Lemma~\ref{lem:infinite}(iv) applies; or $k$ is connected to an even number of vertices, then there is a subdiagram $X'$ which contains $k$ and has the following property: 
the subdiagram $X'$ has at least two cycles and, for any cycle $C$ in $X'$, the product $\prod_{\{i,j\}\in C} (-A_{i,j})$ over all edges of $C$ is positive, so $X'$ is as in Lemma~\ref{lem:infinite}(v) (note that if $k$ is connected to exactly two vertices in $X$, then $X'=Y$), thus $\Gamma$ has an infinite mutation class.

Let us now prove Theorem~\ref{th:acyclic fmc}. If $\Gamma$ is an extended Dynkin (or Dynkin) diagram, then its mutation class is finite by Theorem~\ref{th:mut class ext} and Proposition~\ref{prop:semidefinite comp}(i). For the converse, suppose that $\Gamma$ is a minimal acyclic diagram which is neither Dynkin nor extended Dynkin and let $A$ be an admissible quasi-Cartan companion which is a generalized Cartan matrix. Then $A$ is a generalized Cartan matrix of hyperbolic type \cite[Exercise~4.1]{K}. Thus $A$ is indefinite and non-degenerate \cite[Exercise~4.6]{K}. By Lemmas~\ref{lem:infinite} and ~\ref{lem:C indefinite}, if $\Gamma$ contains an edge whose weight is greater than or equal to $4$ or contains a non-oriented cycle, then it has an infinite mutation class as claimed in the theorem. Let us now assume that $\Gamma$ does not contain any non-oriented cycle and each edge-weight is less than $4$. Then, since $\Gamma$ is not of finite type, there is a sequence of mutations $\mu_k,...,\mu_1$ such that $\Gamma'=\mu_k...\mu_1(\Gamma)$ contains an edge whose weight is at least $4$ or contains  a non-oriented cycle such that for $i=1,...,k-1$, the diagram $\mu_i...\mu_1(\Gamma)$ does not contain any non-oriented cycle nor any edge whose weight is greater than or equal to $4$. Then, by Proposition~\ref{prop:mut adm}, the diagram $\Gamma'$ has an admissible quasi-Cartan companion $A'$ which is mutated from $A$. Since $A'$ is equivalent to $A$, it is non-degenerate. Then, by Lemma~\ref{lem:C indefinite}, the diagram $\Gamma'$, thus $\Gamma$, has an infinite mutation class. This completes the proof of the theorem.

\end{document}